\documentclass[11pt]{amsart}

\usepackage{amsmath,amsthm, amscd, amssymb, amsfonts}

\usepackage{hhline}\usepackage{times} 
\usepackage{graphics}

\usepackage{wasysym}
\input xy
\xyoption{all}
\usepackage{color}

\usepackage{array} 
\usepackage{multirow}

\def\PSL{\operatorname{PSL}}
\def\SL{\operatorname{SL}}

\makeatletter
\def\Ddots{\mathinner{\mkern1mu\raise\p@
\vbox{\kern7\p@\hbox{.}}\mkern2mu
\raise4\p@\hbox{.}\mkern2mu\raise7\p@\hbox{.}\mkern1mu}}
\makeatother
\newcounter{tabla}\stepcounter{tabla}

\newcommand{\Pb}{\boldsymbol{G}}

\def\SO{\operatorname{SO}}
\def\Tt{\operatorname{T}}

\newcommand{\bone}{\boldsymbol{1}}

\newcommand{\Jf}{{\tt{J}}}

\newcommand{\tx}{{\tt{x}}}
\newcommand{\ty}{{\tt{y}}}
\newcommand{\tz}{{\tt{z}}}
\newcommand{\tS}{{\tt{S}}}
\newcommand{\tK}{{\tt{K}}}

\newcommand{\bd}{\mathbf{d}}
\newcommand{\bc}{\mathbf{c}}
\newcommand{\bs}{\mathbf{s}}
\newcommand{\be}{\mathbf{e}}

\newcommand{\Pad}{\overline\G}
\newcommand{\Tad}{\overline\Tt}

\newcommand\sn{\mathbb S_n}

\newcommand{\diag}{\operatorname{diag}}
\newcommand{\Diag}{\operatorname{Diag}}

\newcommand{\Imm}{\operatorname{Im}}

\newcommand\toba{\mathfrak B }

\newcommand{\gr}{\operatorname{gr}}
\newcommand{\trid}{\triangleright}

\newcommand{\Zc}{{\mathcal Z}}

\newcommand{\kk}{\Bbbk}

\newcommand{\Z}{\mathbb Z}
\newcommand{\N}{{\mathbb N}}
\newcommand{\I}{{\mathbb I}}

\newcommand{\G}{{\mathbb G}}
\newcommand{\B}{{\mathbb{B}}}
\newcommand{\T}{{\mathbb{T}}}
\newcommand{\U}{\mathbb{U}}

\newcommand{\F}{{\mathbb F}}

\newcommand{\PGL}{\operatorname{PGL}}
\newcommand{\Fr}{\operatorname{Fr}}

\newcommand{\Oc}{{\mathcal O}}

\newcommand{\Aut}{\operatorname{Aut}}

\newcommand{\Out}{\operatorname{Out}}

\newcommand\card{\operatorname{card}}

\newcommand\Ad{\operatorname{Ad}}

\numberwithin{equation}{section}
\theoremstyle{plain}

\newcommand\id{\operatorname{id}}

\newcommand\Id{\operatorname{Id}}

\def\pf{\begin{proof}}
\def\epf{\end{proof}}

\theoremstyle{definition}
\newtheorem{theorem}{Theorem}[section]
\newtheorem{lemma}[theorem]{Lemma}
\newtheorem{corollary}[theorem]{Corollary}
\newtheorem{proposition}[theorem]{Proposition}
\newtheorem{definition}[theorem]{Definition}

\newtheorem{question}{Question}

\newtheorem*{maintheorem}{Theorem 1.1'}

\theoremstyle{remark}

\newtheorem{remark}[theorem]{Remark}

\def\w{\dot{w}}
\def\v{\dot{v}}

\def\N{\mathbb{N}}
\def\Fr{\operatorname{Fr}}

\def\card{\operatorname{card}}
\def\Z{\mathbb{Z}}
\def\Zc{\mathcal{Z}}

\def\symm{\mathbb{S}}
\def\sym{\mathbb{S}}

\def\F{{\mathbb F}}
\def\Fq{{\mathbb F}_q}
\def\Aut{\operatorname{Aut}}

\def\diag{\operatorname{diag}}
\def\sigmad{\dot{\sigma}}
\def\taud{\dot{\tau}}

\def\int{\mathbb{Z}}

\def\mO{{\mathcal{O}}}

\def\PL{\operatorname{GL}}
\def\Sp{\operatorname{Sp}}
\def\Gal{\operatorname{Gal}}
\def\Tr{\operatorname{Tr}}

\def\gr{\operatorname{gr}}

\newcommand{\ku}{\Bbbk}

\def\Fr{\operatorname{Fr}}

\renewcommand{\eth}{\nu}

\def\pf{\begin{proof}}

\def\epf{\end{proof}}

\def\id{{\rm id}}

\title[$\theta$-semisimple twisted conjugacy classes of type D in $\PSL_n(q)$]{$\theta$-semisimple 
classes of type D in $\PSL_n(q)$}

\author[Carnovale, Garc\'ia Iglesias]{Giovanna Carnovale$^1$, Agust\'in Garc\'ia
Iglesias$^{2}$ }

\email{G.C.: carnoval@math.unipd.it}
\email{A.G.I: aigarcia@famaf.unc.edu.ar}

\thanks{
1 - Dipartimento di Matematica,
Torre Archimede - via Trieste 63 - 35121 Padova - Italy \newline
2 - FaMAF-CIEM (CONICET), Universidad Nacional de C\'ordoba,
Medina A\-llen\-de s/n, Ciudad Universitaria, 5000 C\' ordoba, Rep\'ublica Argentina.
\newline \noindent The work of A.G.I. was partially supported by ANPCyT-FONCyT, CONICET, 
Ministerio de Ciencia y Tecnolog\'ia (C\'ordoba), Secyt (UNC) and the GNSAGA project. 
Part of it was done as a fellow of the Erasmus Mundus EADIC II programme of the EU in the 
Universit\`a degli Studi di Padova. G.C. was partially supported by Progetto di Ateneo CPDA125818/12 and by 
the bilateral agreement between the Universities of C\'ordoba and Padova.}

\date{} 

\begin{document}

\begin{abstract}
Let $p$ be an odd prime, $m\in \N$ and set $q=p^m$, $\Pb=\PSL_n(q)$. Let $\theta$ be a 
standard graph automorphism of $\Pb$, $d$ be a diagonal automorphism and $\Fr_q$ be the 
Frobenius endomorphism of $\PSL_n(\overline{\F_q})$. We show that every $(d\circ \theta)$-conjugacy class 
of a $(d\circ \theta,p)$-regular element in $\Pb$ is represented in some $\Fr_q$-stable maximal torus of  $\PSL_n(\overline{\F_q})$ and 
that most of them are of type D. We write out the possible exceptions and show that, 
in particular, if $n\geq5$ is either odd or a multiple of $4$ and $q>7$, then all such classes are of type D.
We develop general arguments to deal with twisted classes in finite groups.
\end{abstract}

\maketitle
\section{Introduction}

This paper belongs to the series started in \cite{ACGa}, in which we intend to determine
all racks related to (twisted) conjugacy classes in simple groups of Lie type 
which are of type D {\it cf.} \eqref{eqn:typeD}, as proposed in \cite[Question 
1]{AFGaV2}.
This, although being mainly a group-theoretical question, is
intimately related with the classification of finite-dimensional pointed Hopf algebras
over non-abelian groups, see below. In this article we will focus on racks which
arise as non-trivial twisted conjugacy classes in $\PSL_n(q)$ for $q=p^m$, $p$ an odd prime.

\medbreak

Recall that a rack is a non-empty set $X$ together with a binary operation
$\rhd$ satisfying faithfulness and self-distributive axioms, see \ref{subsec:racks}. The
prototypical example of a rack is a twisted conjugacy class $\Oc_x^\psi$ with respect to an automorphism $\psi\in\Aut(G)$ inside a finite group $G$, 
$x\in G$, with
\begin{equation}\label{eqn:twisted-rack}
 y\rhd z=y\psi(zy^{-1}), \qquad y,z\in \Oc_x^\psi.
\end{equation}
This is in fact a quandle, as  $y\rhd y=y$, $\forall\,y\in \Oc_x^\psi$.

A rack $X$ is said to be {\it of type D} when there
exists a decomposable subrack $Y=R\bigsqcup S\subseteq X$and elements $r\in R$, 
$s\in S$ such that $r\rhd(s\rhd (r\rhd s))\neq s$, see Section \ref{subsec:racks}. Their study is deeply connected with the classification problem of finite-dimensional pointed Hopf algebras, as follows.

\medbreak

Let $H$ be a finite dimensional pointed Hopf
algebra over an algebraically closed field
$\ku$ and assume the coradical of $H$ is $\ku G$, for a finite non-abelian group $G$.
Following \cite[Section 6.1]{AG}, there exist a rack $X$ and a 2-cocycle ${\bf q}$ with values in
$\PL(n, \ku)$ such that $\gr H$, the associated graded algebra with respect to the
coradical filtration, contains as a subalgebra the bosonization
$\toba(X,{\bf q})\#\ku G$. See {\it loc. cit.} for unexplained notation. Therefore, it is central for
the classification of such Hopf algebras to know when $\dim\toba(X,{\bf q})<\infty$ for  given
$X$,
${\bf q}$. A rack $X$ is said to {\it collapse} when $\toba(X, {\bf q})$ is infinite
dimensional for any ${\bf q}$. A remarkable result is that {\it if $X$ is of type D, then it collapses}.
This is the content of \cite[Theorem  3.6]{AFGV}, also \cite[Theorem 8.6]{HS1}, 
both of which follow from results in \cite{AHS}.

Now every rack $X$ admits a rack epimorphism $\pi: X \to S$ with $S$ simple and it follows that $X$ is of type 
$D$ if $S$ is so. Hence, determining all simple racks of type D is a drastic reduction indeed for the classification problem, as many groups can be 
discarded and only a few conjugacy classes in simple groups remain. Only for such classes one needs to compute the possible cocycles that yield a finite dimensional Nichols algebras.  Simple racks are classified into three classes \cite{AG}, also \cite{J}, namely {\it affine}, {\it twisted homogeneous} and that of 
{\it non-trivial twisted conjugacy classes on finite simple groups}, see \cite{AG} for definitions. Most (twisted) conjugacy classes in sporadic 
groups are of type D \cite{AFGV2}, \cite{FV}. This is also the case for non-semisimple classes in $\PSL_n(q)$ \cite{ACGa}, for unipotent classes in symplectic groups \cite{ACGa2} and for (twisted) classes in alternating groups \cite{AFGV}. Similar results follow for twisted homogeneous racks \cite{AFGaV}. Affine racks seem to be not of type D.

In this article we begin the analysis of twisted 
classes of type D in $\PSL_n(q)$, for $q$ odd and automorphisms induced by algebraic group  
automorphisms of $\SL_n(\overline{\F_q})$. Recall that the automorphisms in $\PSL_n(q)$ are compositions of  automorphisms induced by 
conjugation in $\PL_n(q)$ (diagonal and inner automorphisms),  powers of a  standard graph automorphism $\theta$ of the Dynkin diagram and powers of the Frobenius 
automorphism $\Fr_p$. Inner automorphisms may be neglected \cite[\S 3.1]{AFGaV}. Diagonal and graph automorphisms are induced by algebraic group  
automorphisms of $\SL_n(\overline{\F_q})$, whereas $\Fr_p$ is induced by an abstract group endomorphism. Their behaviour is therefore different 
\cite[10.13]{steinberg-endo} and 
this is reflected in the structure of the twisted classes. 
In addition, if  the $d\circ \theta^a$-class of $x$ in $\PSL_n(p)$ is of type D, $d$ a diagonal automorphism and $a=0,1$, then 
the $\Fr_p^m\circ\, d\circ \theta^a$-class of $x$ in $\PSL_n(q)$ is 
of type D for every $m$ and every $q$. Thus, we will focus on twisted classes for automorphisms $\psi=d\circ \theta^a$. The analysis of standard conjugacy classes in simple groups of Lie type (corresponding to $a=0$) has been started in \cite{ACGa,ACGa2}. 
For these reasons the first twisted classes to look at are the $\psi$-classes in $\PSL_n(q)$, where $\psi$ is a composition of a diagonal automorphism $d$ 
with $\theta$. In analogy to the case of standard conjugacy classes, it is possible to reduce most of the analysis to the study of classes whose 
behaviour resembles that of semisimple or unipotent ones. However, in contrast 
to that case, the choices to be made depend on the gcd of $|\psi|$ and $p$ {\em 
cf.} Subsection \ref{subsec:semis}. Therefore, the cases of $p$ even and odd must be handled with different methods. The diagonal automorphisms 
always satisfy $(|\psi|,p)=1$ so we restrict to the case $(|\psi|,p)=1$ and we will require  $p$ to be odd.

Set $\Pb=\PSL_n(q)$, $\psi=d\circ\theta\in\Aut(\Pb)$, for $d$ a diagonal automorphism. 
The study of $(\psi,p)$-regular classes in $\Pb$, {\it i.~e.}, of those classes replacing 
semisimple ones, can be reduced to the study of 
$(\theta,p)$-regular $\Pb$-orbits of elements in $\PGL_n(q)$. Such classes 
have a representative in a  maximal torus $\Tad_w^{\Fr_q}$ of 
$\PGL_n(q)$, for some $w\in W^\theta$, where we can take $w$ up to 
conjugation {\it cf.} Theorem \ref{thm:x-in-calT}. It turns out that in most 
cases, the property of being of type D depends on $n$, $q$ and the 
conjugacy class of $w$ in $W^\theta$. Such classes are parametrized by a partition 
$\lambda=(\lambda_1, \dots, \lambda_r)$ of $h=\left[\frac{n}{2}\right]$, with $r\in\N$, $h_i>0$, 
and a certain vector $\varepsilon\in \Z_2^r$. Hence 
our result depends on the number of cycles $r$ of 
$\lambda$ and on the vector $\varepsilon=(\varepsilon_1,\dots,\varepsilon_r)\in 
\Z_2^r$. Let $\bone$ stand for the 
partition 
$(1,\dots,1)$. 
\begin{theorem}\label{thm:one}
Let $q$ be as above. Let $x\in \Tad_w^{\Fr_q}$. Then the class $\Oc_x^{\theta,\Pb}$ is of type D, with the possible exceptions of classes fitting into the 
following table:
\begin{center}
\begin{table}[h]
\begin{tabular}{| c | c |c|  c  | c| m{3cm} |}
  \hline
\multicolumn{3}{|c|}{$w$} & $n$& $q$ & $x$ \\ \hline\hline
\multirow{5}{*}{$\lambda\neq\bone$} & 
$r=2$& $\varepsilon=(0,\varepsilon_2)$
&\mbox{even}&3,5&\mbox{any}\\ \cline{2-6}
& \multirow{4}{*}{$r=1$}&\multirow{2}{*}{$\varepsilon=(0)$}&4&3,7&\mbox{any} \\ 
\cline{4-6}
&&&\raisebox{-.1cm}{4}&\raisebox{-.1cm}{5,9}& 
\raisebox{-.3cm}{$\theta(x)\neq x^{-1}$} \\
\cline{3-6}
&&\multirow{2}{*}{$\varepsilon=(1)$}&\raisebox{-.1cm}{4}&\raisebox{-.1cm}{3,7}&\raisebox{-.2cm}{$\theta(x)\neq x^{-1}$} \\ 
\cline{4-6}
&&&\raisebox{-.1cm}{$2\times$odd}&\raisebox{-.1cm}{any}& 
\raisebox{-.3cm}{$\Oc_x^{\theta,\Pb}\simeq\Oc_{\eth}^{\theta,\Pb}$} \\ 
\hline
\multicolumn{3}{|c|}{\multirow{4}{*}{$\lambda=\bone$}} & any*& 3,5& any* \\ \cline{4-6}
 \multicolumn{3}{|c|}{}&3& 7,13 & any \\ \cline{4-6}
\multicolumn{3}{|c|}{} & 4& $\equiv3(4)$ & any \\ \cline{4-6}
\multicolumn{3}{|c|}{} & 4& 9 & any \\ \hline
\end{tabular}
\medskip
\caption{Possible exceptions; $\nu$ as in \eqref{eqn:missing}.}\label{table:1}
\end{table}
\end{center}
\end{theorem}
* Actually, some of the classes listed on the table are of type D, for instance 
when $n\geq 6$, $n\neq 7$ and $\varepsilon=(0,\ldots,0)$, see Lemma \ref{lem:weyl-j}. See 
also Remark \ref{rem:r=2}.

\medskip

We present this result in the language of Nichols algebras, as a partial answer 
in this cases to \cite[Question 2]{AFGaV2}, see also
\cite[Theorem 3.6]{AFGV}, and {\it loc. cit.} for unexplained notation. 
Consider the classes $\Oc_x^{\theta,\Pb}$ 
in Theorem \ref{thm:one} as racks with the rack structure 
\eqref{eqn:twisted-rack}. These are simple racks.
\begin{corollary}
Let $X=\Oc_x^{\theta,\Pb}$, $x\in \Tad_w^{\Fr_q}$. Then $\dim 
\toba(\Oc_x^{\theta,\Pb},{\bf q})=\infty$ for any cocycle ${\bf q}$ on $X$, with the 
possible exceptions of the classes in Table \ref{table:1}.\qed
\end{corollary}

Also, an extract of Theorem \ref{thm:one} can be rephrased as follows.
\begin{maintheorem}
Let $p$ be an odd prime, $m\in\N$, $q=p^m$. Set $\Pb=\PSL_n(q)$, $\psi=d\circ\theta\in\Aut(\Pb)$, for $d$ a diagonal automorphism. 

If $n\geq 5$, $q\geq 7$, then any $(\psi,p)$-regular class $\Oc$ is of type D 
with the possible exception $n=2\times$odd, $\Oc\simeq 
\Oc_1^{\Ad(\eth^{-1})\circ\theta,\Pb}$, 
$\eth$ as in \eqref{eqn:missing}.\qed
\end{maintheorem}

When $\psi=\theta$, we obtain the following for classes with trivial 
$(\theta,p)$-regular part (also called $\theta$-semisimple part) which is the content of Propositions 
\ref{pro:trivial-even} and \ref{pro:trivial-odd}:
\begin{proposition}
Let $\Oc$ be a $\theta$-twisted conjugacy class with trivial $\theta$-semisimple part. Then $\Oc$ is of type D provided 
\begin{enumerate}
\item $n>2$ is even, the unipotent part is nontrivial, and $q>3$.  
\item $n>3$ is odd and the Jordan form of its $p$-part in $\Pb^\theta$ corresponds to the partition $(n)$.\qed
\end{enumerate}
\end{proposition}

The paper is organized as follows. In Section \ref{sec:preliminaries} we fix the notation and recall some generalities about racks and the group $\PSL_n(q)$. In Section 
\ref{sec:general} we discuss some general techniques to deal with twisted conjugacy classes in a finite group. In Section \ref{sec:PSL} we focus on $\PSL_n(q)$ and we begin a systematic approach to the study of its twisted classes, that includes an analysis of the Weyl group. In Section \ref{sec:regular} we concentrate on $\theta$-semisimple classes and obtain the main results of the article. In Section \ref{sec:unipotent} we present some results on classes with trivial $\theta$-semisimple part.

\section{Preliminaries}\label{sec:preliminaries}

Let $H$ be a group, $\psi\in \Aut(H)$. A $\psi$-twisted conjugacy class, 
or simply a twisted conjugacy class, is an orbit for the action of $H$ on itself 
by $h\cdot_\psi x=hx\psi(h)^{-1}$. We denote this class by $\mO_h^\psi$.  
If $K< H$ is $\psi$-stable, we will write $\mO_h^{\psi,K}$ to denote the orbit 
of $h$ under the restriction of the $\cdot_\psi$-action to $K$. In particular, 
$\mO_h=\mO_h^\id$ denotes the (standard) conjugacy class of $h\in H$. 
The stabilizer  in $K<H$ of an element $x\in H$ for the twisted action will be 
denoted by $K_\psi(x)$
so that $H_{\id}(x)$ is $H_x$, the usual centralizer of $x$. For any 
automorphism
$\psi$ of a group $H$, we write $H^\psi$ for the set of $\psi$-invariants 
in $H$.  The inner automorphism
given by conjugation by $x\in H$ will be denoted by $\Ad(x)$. If 
$K\triangleleft H$ is normal, then we also denote by $\Ad(x)$ the automorphism 
induced from the conjugation by $x\in H$. 
$\Zc(H)$ will denote the center of $H$. Recall that the group $\mu_n(\F_q)$ of roots of unity in a finite field $\F_q$ is isomorphic to  $\Z_{\bd}$,
for $\bd:=(n, q-1)$.

We denote by 
$\sn$, $n\in\N$, the symmetric group on $n$ letters.
We also set  $\I_n:=\{1,2,\,\ldots,\,n\}$ and 
$(b)_a=1+a+a^2+\dots+a^{b-1}$, $a,b\in \N$.

\subsection{Racks}\label{subsec:racks}

A \emph{rack} $(X, \rhd)$ is a non-empty finite set $X$
together with a function $\rhd:X\times X\to X$ such that $i\rhd
(\cdot):X\to X$ is a bijection for all $i\in X$ and $$i\rhd(j\rhd
k)=(i\rhd j)\rhd (i\rhd k), \, \forall i,j,k\in X.$$
Recall that a rack $(X, \rhd)$ is a \emph{quandle} when $i\trid i=i$, $\forall\,i\in X$.

We shall write simply $X$ when the function $\rhd$ is clear from the context. 

\medbreak

If $H$ is a group, then the conjugacy class $\mO_h$ of any element $h\in H$ is a 
rack, with the function $\rhd$ given by conjugation. More generally,  if 
$\psi\in\Aut(H)$, any twisted conjugacy class in $H$ is a rack with rack 
structure given by \eqref{eqn:twisted-rack},
see \cite[Theorem 3.12, (3.4)]{AG}. These are indeed examples of quandles. 

\medbreak

A subrack $Y$ of a rack $X$ is a subset $Y\subseteq X$ such that $Y\rhd 
Y\subseteq Y$. A rack is said 
to
be \emph{indecomposable} if it cannot be decomposed as the disjoint union of two 
subracks. A rack $X$ is said to be \emph{simple} if $\card X > 1$ and for any 
surjective 
morphism of racks $\pi: X \to Y$, either $\pi$ is a bijection or $\card Y = 1$.

\subsubsection{Racks of type D}

A rack $X$ is of type D when there
exists a decomposable subrack $Y=R\bigsqcup S$ of
$X$ and elements $r\in R$, $s\in S$ such that
\begin{equation}\label{eqn:typeD}
r\rhd(s\rhd (r\rhd s))\neq s.
\end{equation}

If a rack  $X$
has a subrack of type D, or if there is a rack epimorphism 
$X\twoheadrightarrow Z$ and $Z$ is of type D, then $X$ is again so. In 
particular, if $X$ is decomposable and $X$ has a component of type D, then 
$X$ is of type D. On the other hand, if $X$ is indecomposable, then it admits a projection 
$X\twoheadrightarrow Z$, with $Z$ simple. Hence, in the quest of racks of 
type D it is enough to focus on simple racks. The classification of simple 
racks is given in \cite[Theorems 3.9, 3.12]{AG}, see also
\cite{J}. A big class consists of twisted conjugacy classes in finite simple 
groups.

\begin{remark}\label{rem:typeD}
Let $\Oc$ be a $\psi$-twisted conjugacy class. Then $\Oc$ is of type D if there are $r,s\in \Oc$ 
such that $r\notin \Oc_s^{\psi,L}$, for $L$ the $\psi$-stable closure of 
the subgroup generated by $r$ and $s$, and
\begin{equation}\label{eq:rs}
 r\psi(s)\psi^2(r)\psi^3(s)\neq s\psi(r)\psi^2(s)\psi^3(r).
\end{equation}
In fact, if the above conditions hold, we set $S=\mO^{\psi,L} _s$ and $R=\mO^{\psi,L} 
_r$ and then $Y=R\bigsqcup S$ is a decomposable subrack of $\Oc$ which satisfies 
\eqref{eqn:typeD}.

If $\psi=\id$ then the condition is also necessary: if $\Oc$ is of type D, then 
there are $r,s\in\Oc$, $r\notin \Oc_s^{L}$, satisfying \eqref{eq:rs} 
\cite[Remark 2.3]{ACGa}.
\end{remark}

\subsection{The group $\Pb=\PSL_n(q)$}\label{sec:Gb}
Fix $n\in\N$. Let $p\in\N$ be a prime number and let $\kk=\overline{\F_p}$. Fix $m\in\N$, $q=p^m$. We assume throughout the paper that $n>2$ or $q\neq 2,3$.

We fix once and for all the following notation:
\begin{align}\label{eqn:gps}
 &\G=\SL_n(\kk), &&\Pad=\PSL_n(\kk), && \Pb:=\PSL_n(q).
\end{align}
We also fix $\pi\colon  \PL_n(\kk)\to \PGL_n(\kk)\simeq \Pad$ the usual projection. We shall keep the name $\pi:=\pi_{|\G}\colon \G \to \Pad$ for the restriction of $\pi$ to $\G$. 
We fix the subgroups of diagonal matrices
\begin{align}\label{eqn:tori}
&\Tt\leq \PL_n(\kk), &&\T\leq\G, &&\Tad:=\pi(\Tt)\leq \Pad.
\end{align}
\subsubsection{General properties of $\Pb$}
Consider  the exact sequence:
\begin{align}\label{eqn:exact}
 & 1\longrightarrow \Zc(\G)\longrightarrow \G\overset{\pi}\longrightarrow 
\Pad\longrightarrow 1
\end{align}
and let $F=\Fr_p^m$ be the endomorphism of $\PL_n(\kk)$ raising every entry in 
$X\in\PL_n(\kk)$ to the $q$-th power.  Taking $F$-points, \eqref{eqn:exact} 
yields:
\begin{align*}
&1\longrightarrow \Zc(\SL_n(q))\longrightarrow \SL_n(q)\longrightarrow \PGL_n(q).
\end{align*}
Then $\Pb\leq \PGL_n(q)$ is the image of the last arrow: $$\Pb=\PSL_n(q)\simeq 
\SL_n(q)/\Zc(\SL_n(q))\simeq \SL_n(q)/\Z_{\bd},$$ for $\bd=(n,q-1)$.
The group $\Pb$ is simple\footnote{Recall that $\PSL_2(2)\simeq 
\mathbb{S}_3$, $\PSL_2(3)\simeq \mathbb{A}_4\leq\mathbb{S}_4$.}.

\smallbreak
We will denote by $\B, \U, \U^-\leq\G$ be the subgroups of $\G$ of upper 
triangular, unipotent upper-triangular, unipotent lower-triangular matrices. Set 
 $$W:=N_{\G}(\T)/\T\simeq N_{\Pad}(\Tad)/\Tad \simeq 
\sn.$$

Recall that $[\SL_n(q),\SL_n(q)]=\SL_n(q)$ and 
$[\PGL_n(q),\PGL_n(q)]=\Pb$, for $n>2$ or $q\neq 2,3$. Also, we have the 
identifications:
\begin{multline*}
\Pad^F=\PGL_n(q)=\Tad^F [\PGL_n(q),\PGL_n(q)]=\Tad^F \Pb\\
\simeq 
\PL_n(q)/\Zc(\PL_n(q))\simeq \PL_n(q)/\F_q^\times.
\end{multline*}

\subsubsection{Automorphisms of $\Pb$}

Recall that a {\it diagonal automorphism} of $\Pb$ is an automorphism induced 
by conjugation by an element in $\Tad^F$. The {\it graph automorphism} 
$\theta\colon \PL_n(\kk)\to \PL_n(\kk)$ is given by 
$x\mapsto \Jf_n\,^t x^{-1}\Jf_n^{-1}$, for  
\begin{align}\label{eqn:J}
 \Jf_n=\left(\begin{smallmatrix}
0&\dots&0&1\\ 0&\dots &-1&0\\
\vdots&\Ddots&\vdots&\vdots\\
(-1)^{n-1}&\dots&0&0\end{smallmatrix}\right).
\end{align}
It induces a non-trivial automorphism of $\G$ for $n\geq3$ and it is unique up 
to inner automorphisms\footnote{Indeed, this 
is not the choice made in \cite{ACGa} but it is, however, more adequate for our 
setting.}. It also induces automorphisms of $\PL_n(q), \SL_n(q)$, $\PGL_n(q)$ and $\Pb$. We will drop 
the subscript $n$ and write $\Jf=\Jf_n$ when it can be deduced from the context.

\medbreak

By \cite[Theorem 24.24]{MT} every automorphism of $\Pb$ is the composition of an inner, 
a diagonal, a power of $\Fr_p$ and a power of $\theta$, so the elements in 
group of outer automorphisms of $\Pb$ have representatives  in
$\Out(\Pb):=\langle \Fr_p,\theta, \Ad(t)\,:t\in \Tad^F\rangle$.

\section{General arguments}\label{sec:general}

In this section we present some general techniques to deal with twisted 
conjugacy classes in finite groups. 

\smallbreak

We start with a well-known lemma.

\begin{lemma}
\label{lem:inner} 
Let $H$ be a finite group, $\varphi\in \Aut(H)$. Let  $K,N< 
H$ be $\varphi$-stable subgroups, with $N\lhd H$. Fix $x\in H$.

(1) The set $\Oc_x^{\varphi, K}$ is a subrack of $\Oc_x^{\varphi, H}$ if and 
only if for every $k\in K$ there is $t\in H_\varphi(x)$ such that $xkx^{-1}t\in K$.

(2)\cite[\S 3.1]{AFGaV} Assume $\varphi=\Ad(x)\circ\psi$, for some $\psi\in 
\Aut(H)$. Then for every $g\in H$ there are racks isomorphisms
$\Oc_g^{\varphi,H}\simeq \Oc_{gx}^{\psi,H}$ and 
$\Oc_g^{\varphi,N}\simeq\Oc_{gx}^{\psi,N}$. 

(3) Let $y\in H$ with $y\in\Oc_x^{\varphi,H}$. Then $\Oc_x^{\varphi,N}\simeq
\Oc_y^{\varphi,N}$.
\end{lemma}
\pf
(1) is straightforward. In (2), we have the equality of sets 
$\Oc_g^{\varphi,H}=\Oc_{gx}^{\psi,H}x^{-1}$ and right multiplication by $x$ defines the 
rack isomorphism. The second isomorphism follows by 
restriction. As for (3), let $g\in H$ be such that $g\cdot_\varphi x=y$. Then the 
map $T:\Oc_x^{\varphi,N}\to \Oc_y^{\varphi,N}$ given by $T(z)=g\cdot_\varphi z$ is a rack 
isomorphism. Observe that if $z=h\cdot_\varphi x$ then $T(z)=(ghg^{-1})\cdot_\varphi y$.
\epf  
\begin{remark}\label{rem:comm}
Notice that the assumption in (1) in Lemma \ref{lem:inner} holds if $x\in 
N_H(K)$. In particular, it always holds for $K\lhd H$. Also, (2) allows us to 
neglect inner automorphisms of $H$. 
\end{remark}

The following slight generalization of \cite[Lemma 2.5]{FV} will be very useful.
\begin{lemma}\label{lem:vendra}
Let $H$ be a finite group and let $K\lhd H$. Let $s\in H$ be a non-trivial  
involution. Then $\Oc_s^K$ is a rack of type D if and only if there is $r$ in 
$\Oc_s^K$ such that $|rs|$ is even and greater than $4$.
\end{lemma}

\pf By Lemma \ref{lem:inner}, Remark \ref{rem:comm}, $\Oc_s^K$ is a rack.
Observe first that, if $r\in \Oc_s^K$, then the racks
$\Oc_s^{\langle r,s\rangle}$ and $\Oc_r^{\langle 
r,s\rangle}$ are subracks of $\Oc_s^K$. Indeed, if $r=k\trid s=ksk^{-1}$, then 
a generic element of $\langle s,\, r\rangle$ has the form
$y_{a,b}=s^a ksk^{-1}s\cdots ksk^{-1}s^b$ for $a,b\in\{0,1\}$.
Let $sks=l\in K$. Then, if $a=1$ we have
\begin{align*}
&y_{1,b}\trid s=y_{1,0}\trid s=lk^{-1}\cdots l k^{-1}\trid s\in \Oc_s^K,\\
&y_{1,b}\trid r=lk^{-1}\cdots l k^{-1} s^b k s^b\trid s\in\Oc_s^{K},
\end{align*}
whereas if $a=0$ we have 
\begin{align*}
&y_{0,b}\trid s=y_{0,1}\trid s=kl^{-1}\cdots k l^{-1}\trid s\in \Oc_s^K,\\
&y_{0,b}\trid r=kl^{-1}\cdots kl^{-1}s^{b-1}ks^{b-1}\trid s\in \Oc_s^K,
\end{align*} so
the racks $\Oc_s^{\langle r,s\rangle}, \Oc_r^{\langle r,s\rangle}\subset \Oc_s^K$.
Now, if an $r$ as in the statement exists, then $r\trid(s\trid(r\trid s))\neq s$ and
$\Oc_s^{\langle r,s\rangle}$ and $\Oc_r^{\langle r,s\rangle}$ are disjoint,  so 
$\Oc_s^K$ is of type D by Remark \ref{rem:typeD} for $\psi=\id$.
Conversely, if there is no such an $r$, then for every $x\in \Oc_s^K$ either 
$|xs|\leq 4$ or
it  is odd, so either $(xs)^2=(sx)^2$ or $\Oc_s^{\langle 
s,x\rangle}=\Oc_x^{\langle s,x\rangle}$ and Remark \ref{rem:typeD} for 
$\psi=\id$ applies once more.
\epf

\begin{remark}\label{rem:reductions}Let $H$ be a finite group,  $\phi\in\Aut(H)$, $h\in 
H$.
\begin{enumerate}
\item Assume $K=H_h$ is $\phi$-stable. If $k\in K$, then 
$\Oc_{kh}^{\phi,K}= \Oc_k^{\phi, K} h$ as sets  and right multiplication by $h^{-1}$ gives a 
rack isomorphism $\Oc_{kh}^{\phi,K}\simeq\Oc_k^{\phi,K}$.

\item Let  $L=H\rtimes\langle \phi\rangle$. Then, for $x=g\phi$,  we have the equality of 
sets: $\Oc_g^{\phi,H}=\Oc_x^L\,\phi^{-1}$ and $y\mapsto y\phi$ induces a rack 
isomorphism $\Oc_g^{\phi,H}\simeq\Oc_x^L$.
\end{enumerate}
\end{remark}

\begin{remark}\label{rem:reductions2} Let $H$ be a finite group,  $\phi\in\Aut(H)$. Let 
$A$ be a $\phi$-stable abelian subgroup of $H$, $a\in A$.
\begin{enumerate}

\item\label{item:abelian} By Remark \ref{rem:reductions} (1), $\Oc_a^{\phi, A}\simeq 
\Oc_1^{\phi,A}$ as racks. Moreover $\gamma\colon A\to A$, $b\mapsto 
b\phi(b^{-1})$, is a group morphism and 
$\Oc_1^{\phi,A}=\Imm(\gamma)\simeq A/A^\phi$ as groups. 

\item\label{item:abelian2} If  
$\phi$ is an involution, then $\Oc_a^{\phi, A}$ is of type D if and only if 
there is $b\in 
A/A^\phi$ such that $|b|$ is even, $>4$ by Remark \ref{rem:reductions} (2) and  
Lemma \ref{lem:vendra}.

\item\label{item:p-part} Let $p$ be a prime number dividing $|H|$. Let 
$h=us=su\in H$ be the (unique) 
decomposition of $h$ as a product of a $p$-element $u$ and a $p$-regular 
element $s$. If  $\Oc_{u}^{H_s}$  is of type D, then $\Oc_h$ is again so, as  
$\Oc_{u}^{H_s}$
identifies with a subrack of $\Oc_h^H$.
 \end{enumerate}

\end{remark}

\begin{remark}\label{rem:stable-orbit}
Let $H$ be a group, let $\phi,\psi\in \Aut(H)$, with $\phi\psi=\psi\phi$, and 
let $N\lhd H$ be $\phi$-stable and $\psi$-stable.

(1) If $\Oc_{h}^{\phi,N}\cap H^\psi\neq\emptyset$, then 
$\psi(\Oc_{h}^{\phi,N})=\Oc_{h}^{\phi,N}$. 
Indeed, let $x\in \Oc_t^{\phi,N}$ with 
$\psi(x)=x$. Now, if $y=kx\phi(k^{-1})\in 
\Oc_x^{\phi,N}=\Oc_h^{\phi,N}$, $k\in N$, then 
$\psi(y)=\psi(k)x\phi(\psi(h)^{-1})\in \Oc_h^{\phi,N}$.

(2) Conversely, if $\psi(\Oc_{h}^{\phi,N})=\Oc_{h}^{\phi,N}$ and the map $N\to N$, given by $x\mapsto x^{-1}\psi(x)$, $x\in N$,
is surjective, then $\Oc_{h}^{\phi,N}\cap H^\psi\neq\emptyset$. 
To see this, fix $g\in N$ such that $\psi(h)=gh\phi(g^{-1})$ and let $x\in N$ be such that $g^{-1}=x^{-1}\psi(x)$. 
Then it follows that $x\cdot_\phi h\in H^\psi\cap \Oc_{h}^{\phi,N}$.

\subsection{$(\psi,p)$-elements and $(\psi,p)$-regular elements}\label{subsec:semis}

Let $H$ be a finite group,
$p$ be a prime number dividing $|H|$ and let $\psi\in \Aut(H)$, with $\ell:=|\psi|$. Set 
$\widehat{H}=H\rtimes\langle\psi\rangle$.

\begin{definition}
An element $h\in H$ is called $(\psi,p)$-regular if $h\psi$ is $p$-regular in $\widehat{H}$, {\it i.~e.}
if $(|h\psi|,p)=1$.
An element $h\in H$ is called a $(\psi,p)$-element if $h\psi$ is a $p$-element in $\widehat{H}$, {\it i.~e.}  if $|h\psi|=p^a$ for some $a\in\N$. \end{definition}

Let $\psi=\psi_r\psi_p$ be the decomposition of $\psi$ as a product of its usual $p$-regular part and its $p$-part in 
$\Aut(H)$. 
Then for every $h\psi$  in $\widehat{H}$ we have $h\psi=s\psi_r(u)\psi=u\psi_p(s)\psi$ where $s$ is $(\psi_r,p)$-regular 
 and $u$ is a $(\psi_p,p)$-element in  $H$.

In the quest of $\psi$-classes of type D, a first analysis  can be done by looking at subracks given by the orbits with 
respect to  $H^{\psi_r}$ or $H^{\psi_p}$.
For this reason, the analysis should begin with the cases in which either $\psi_p=1$, {\it i.~e.} when $(\ell,p)=1$, 
or when  $\psi_r=1$, {\it i.~e.} when $\ell$  is a power of $p$.

If $(\ell,p)=1$, then for every $h\in H$ there is a unique decomposition $h=us=s\psi(u)$ with $u$ a $p$-element in $H$
and $s$ a $(\psi,p)$-regular element. In this case $s$ is $(\psi,p)$-regular if and only if the norm ${\rm Norm}_\psi(s):= 
s\psi(s)\cdots\psi^{\ell-1}(s)$ is $p$-regular in $H$. Here, if $C=H_{\psi}(s)$ and $C'=\widehat{H}_{s\psi}$, then Remarks \ref{rem:reductions} (2) and \ref{rem:reductions2} \eqref{item:p-part} give the rack inclusions 
\begin{equation}\label{eqn:unipotent-inclusions}
\Oc_h^{\psi,H}\simeq \Oc_{h\psi}^{\widehat{H}}\supset  \Oc_{u}^{C'}\supset\Oc_u^{C}. 
\end{equation}
So if  $\Oc_{u}^{C}$  is of type D, then $\Oc_h^{\psi,H}$ is again so. Hence the 
first classes to be attacked are either standard conjugacy classes of 
$p$-elements in $C$ or twisted classes of $(\psi,p)$-regular elements in $H$. 
The latter are dealt with in Section \ref{sec:regular}.

Similarly, if $\ell=p^b$ for some $b>0$, then for each $h\in H$ there is a 
unique decomposition $h=su=u\psi(s)$ with $s$ a usual $p$-regular element in $H$ and $u$ 
a $(\psi,p)$-element. In this case $u$ is a $(\psi,p)$-element if and only if 
${\rm Norm}_\psi(u)$ is a $p$-element in $H$. The first reduction 
is to look at classes of $(\psi,p)$-elements and the standard $p$-regular 
classes in $H_{\psi}(u)$.  We will not pursue this analysis in this paper.

\smallbreak

Notice that, when dealing with twisted classes in simple groups of Lie type, there is a privileged choice for $p$, namely, the 
defining characteristic.

\section{Twisted classes and $\PSL_n(q)$}\label{sec:PSL}

In this section we collect some results  that  contribute to establish a systematic approach to twisted classes in 
$\PSL_n(q)$. This in particular requires a detailed study of the conjugacy classes in the subgroup of 
$\theta$-invariant elements of the Weyl group, and of the corresponding $F$-stable maximal tori in $\Pb$, that we develop in \S \ref{subsec:weyl}.

Recall the notation from \S \ref{sec:Gb}, specially in \eqref{eqn:gps}, \eqref{eqn:tori}.
Next proposition deals with diagonal automorphisms $d=\Ad(t)$, $t\in \Tad^F$.

\begin{proposition}\label{pro:x-in-PGL}
Let $x\in\Pb$, $\varphi=\Ad(t)\circ \psi\in\Aut(\Pb)$, $t\in \Tad^F$. Let $y=t^{-1}x\in\Pad^F$.  Then  
$\Oc_x^{\varphi,\Pb}\simeq \Oc_y^{\psi,\Pb}$. If, in addition,  $\psi\in\Out(\Pb)$ and $z\in 
\Oc_y^{\psi,\PGL_n(q)}$, then $\Oc_x^{\varphi,\Pb}\simeq \Oc_z^{\psi,\Pb}$.
\end{proposition}

\pf
In this case, $x=ty$ and $\Oc_x^{\varphi,\Pb}\simeq \Oc_y^{\psi,\Pb}$ by Lemma 
\ref{lem:inner} (2). The last assertion is Lemma \ref{lem:inner} (3).
\epf

Let $\psi=\Fr_p^a\circ \theta^b\in \Aut(\PL_n(q))$ and let  $\ell:=|\psi|$. 
Then $\psi$ induces an automorphism of $\SL_n(q), \PSL_n(q)$ and $\PGL_n(q)$ of the same order. 
Let $H$ be either $\PL_n(q)$, $\SL_n(q), \PSL_n(q)$, or $\PGL_n(q)$, $\widehat{H}=H\rtimes\langle \psi\rangle$.

If $(\ell,p)=1$, then the $(\psi,p)$-elements in $H$ are the unipotent elements in $H$. The $(\psi,p)$-regular elements 
are those $g\in H$ such that ${\rm Norm}_\psi(g)$ is semisimple.
If, instead, $\ell=p^b$ for some $b>0$, then the $(\psi,p)$-regular elements  in $H$ are the semisimple elements in $H$, 
 while the $(\psi,p)$-elements are those $g\in H$ such that 
${\rm Norm}_\psi(g)$ is a $p$-element.

We will concentrate on the case $(\ell,p)=1$. We have the following 
equivalence.

\begin{lemma}\label{lem:equivalence-1}
Let $\psi\in \Aut(\PL_n(q))$ with $(|\psi|,p)=1$.
Then $\tx\in \PL_n(q)$ is $(\psi,p)$-regular if and only if $x=\pi(\tx)\in \PGL_n(q)$ 
is $(\psi,p)$-regular.
 \end{lemma}
\pf
${\rm Norm}_\psi(x)$ is semisimple if and only $\pi({\rm Norm}_\psi(x))={\rm Norm}_\psi(\tx)$ is so.
\epf

\subsection{The case $\psi=\theta$, $p\neq 2$}

We intend to study twisted classes for automorphisms induced from algebraic group automorphisms. 
By Remark \ref{rem:comm} and Proposition \ref{pro:x-in-PGL}, we may reduce to the case $\psi=\theta$.
We will focus on the case of $p$ odd and we shall investigate $(\psi,p)$-regular classes.
 
 \begin{remark}\label{rem:vinberg}It was pointed to us by Prof. Vinberg that when the group 
 is $\PL_n({\overline{F_q}})$ and $\psi=\theta$, then the map 
 $x\mapsto x\Jf$ allows to identify the $\theta$-twisted 
 conjugacy class of $x$ with the equivalence classes of the non-degenerate bilinear form with 
 associated matrix $x\Jf$. Thus, the classification of twisted classes in this case can be 
 deduced from the classification of bilinear forms on $\overline{\F_q}^n$. The latter, in turn, goes over 
 in odd characteristic, as the classification in characteristic zero which is to be found for instance 
 in \cite{HoP}. From this, $\SL_n(q)$-orbits could be also classified. However,
 since the action of the center by twisted conjugation is non-trivial, the step to $\PSL_n(q)$-orbits of 
 elements in $\PGL_n(q)$ would need slight care. The main reason for our apparently less natural approach is 
 related to the 
 general problem of detecting twisted classes of type D in all finite simple groups.  
 One of the aims in this paper is to propose  a general systematic approach that could be applied, 
 at least,  to all finite simple groups of Lie type. 
\end{remark}

\begin{lemma}\label{lem:equivalence}
Let $\tx\in \PL_n(q)$.
\begin{enumerate}
\item $\tx$ is $\theta$-semisimple if 
and only if there is a $g\in\PL_n(\kk)$ such that $g\cdot_\theta \tx$ lies in a 
$\theta$-stable torus $\Tt_0$ in $\PL_n(\kk)$.
\item $\tx$ is $\theta$-semisimple if and only if there is a $g'\in\SL_n(\kk)\subset \PL_n(\kk)$ such that 
$g'\cdot_\theta \tx\in \Tt$.
\end{enumerate}
 \end{lemma}
\pf
(1) is \cite[Proposition 3.4]{mohr}. Following the construction in \cite[page 382]{mohr} 
we can make sure that $\Tt_0$ is $F$-stable and that it is contained in $\Tt$. For (2), 
let $\Zc:=\Zc(\PL_n(\kk))$, hence $\PL_n(\kk)=\Zc\G$ and $\theta$ acts
as inversion on $\Zc$. Therefore, if $z\in \Zc$, then $z\cdot_\theta x=xz^2$.
Let $g=zg'\in \Zc\,\G$ be such that $g\cdot_\theta x=t\in \Tt$. Then $g'\cdot x= t z^{-2}\in \Tt$, 
as $\Zc$ is contained in every maximal torus.
\epf

The lemma above motivates the following definition.

\begin{definition}
We say that an element $x\in\PGL_n(q)$ is $\theta$-semisimple if it is $(\theta,p)$-regular.
\end{definition}

\end{remark}

\subsection{$F$-stable maximal tori}\label{subsec:weyl}

In this section we collect preparatory material in order to find suitable representatives 
of $\Pb$-classes of  $\theta$-semisimple elements in $\PGL_n(q)$. Unless otherwise stated, $p$ is arbitrary. 

Let $H$ denote either $\G$, $\Pad$ or $\PL_n(\kk)$ and, consequently, set  
$K=\T$, $\Tad$ or $\Tt$ ($=\T\Zc(\PL_n(\kk))$).
Let $w \in W$, $\w\in wK$ and $g=g_w\in H$ 
be such that $g^{-1}F(g)=\w$ (Lang-Steinberg's Theorem). We set 
\begin{equation}
K_w:=gKg^{-1}. 
\end{equation}
Then $K_w$ is an $F$-stable maximal torus of $H$ and all $F$-stable maximal tori in 
$H$ are obtained this way \cite[Proposition 25.1]{MT}. Two tori $K_w$ and $K_\sigma$ are $H^F$-conjugate if and only 
if $\sigma$ and $w$ are $W$-conjugate.  We will provide a $\theta$-invariant version of this 
fact in Lemma \ref{lem:invt} for $K=\Tt$ and $\Tad$. 
We set
\begin{align}\label{eq:invariance}
 &F_w:=\Ad(\w)\circ F, \text{ so } (K_w)^F=g K^{F_w}g^{-1}.
 \end{align}
The automorphisms $\theta$ and $F$ preserve $\T$, hence they induce automorphisms on $W$ 
which we denote by the same symbol.  
The action of $F$ on $W$ is trivial, whereas the action of $\theta$ is conjugation by the 
longest element $w_0\in W$, so $W^\theta=W_{w_0}$. Observe that
\begin{align*}
w_0=\begin{cases}(1,\,n)(2\,n-1)\dots(h,\,h+1)&\mbox{ if }n=2h,\\
(1,\,n)(2,\,n-1)\dots(h,\,h+2)&\mbox{ if }n=2h+1.
\end{cases}
\end{align*}
Any $\sigma\in W^\theta$  can be written as $\sigma=\omega\tau$ where $\omega$ 
permutes the $2$-cycles in $w_0$ and $\tau$ is a product of  transpositions occurring in the cyclic decomposition of $w_0$. 
In fact, $W^\theta\simeq \sym_h\rtimes\Z_2^h$, where 
$h=\left[\frac{n}{2}\right]$, the elements in  $\sym_h$ correspond to products 
$c \theta(c)$ where $c$ is a cycle in $\sym_{\I_h}\leq\sym_n$, $\theta(c)=w_0 c 
w_0$ and the elements in $\Z_2^h$ are products of transpositions of the form 
$(i, n+1-i)$. 

\begin{remark}\label{rem:representatives}
There is a set of representatives $\{\sigmad\}\subset N_{\G}(\T)$ of $W$ such 
that $\sigmad\in N_{\G}(\T)^\theta$ if $\sigma\in W^\theta$, \cite[8.2, 8.3 (b)]{steinberg-endo}. In addition,
$\G^\theta=\Sp_n(\kk)$ if $n$ is even, $\G^\theta=\SO_n(\kk)$ if $n$ is odd and $W^\theta$ is the corresponding Weyl 
group.
\end{remark}

\begin{lemma}\label{lem:invt}
Let $w,\sigma\in W^\theta$. Then $\Tt_w$ and $\Tt_\sigma$ are  $\SL_n(q)^\theta$-conjugate 
if and only if $\sigma\in\Oc_w^{W^\theta}$ if and only if $\overline{\Tt}_w$ and $\overline{\Tt}_\sigma$ are  $\pi(\SL_n(q)^\theta)$-conjugate.
\end{lemma}
\pf Since ${\rm Ker}(\pi)$ consists of central elements, it is enough to prove the first equivalence. 
By Remark \ref{rem:representatives} there are representatives $\w$, $\sigmad$ 
 of $w$ and $\sigma$ in $\G^\theta\cap N(\Tt)$. By Lang-Steinberg's Theorem applied to $\G^\theta$ we may find 
 $y,z\in\G^\theta$ such that $y^{-1}F(y)=\w$,  $z^{-1}F(z)=\sigmad$.
 
Assume there is $x\in \SL_n(q)^\theta$ such that $x\Tt_wx^{-1}=\Tt_\sigma$. Then, 
$\taud:=z^{-1}xy\in N(\Tt)\cap\G^\theta$ and $\taud\w \taud^{-1}\Tt=\taud\w F(\taud^{-1})\Tt=\sigmad \Tt$. 

Conversely, assume there is $\tau\in W^\theta$ such that $\tau w\tau^{-1}=\sigma$. Let 
$\taud\in \G^\theta\cap \tau \Tt$. Then there exist $h,k\in \G^\theta\cap 
\Tt=\Tt^\theta=\Tt^{\theta,\circ}$ such that $F(\taud)=\taud h$ and  $\sigmad=\taud 
\w\taud^{-1}k$. For $t\in \Tt^{\theta}$ we set $x_t=z\taud t y^{-1}\in \G^\theta$. 
Now, 
$x_t \Tt_w x_t^{-1}=\Tt_\sigma$. In addition, $x_t\in \SL_n(q)^\theta$ if and only if $t=\w(\taud^{-1}k\taud)h F(t)\w^{-1}$. This happens if 
and only if  $t^{-1} (\Ad(\w)\circ F)(t)=\w h^{-1}(\taud^{-1}k^{-1}\taud)\w^{-1}$. By Lang-Steinberg's 
Theorem applied to the Steinberg  endomorphism $\Ad(\w)\circ F$ on $\Tt^{\theta}$, there is 
$t\in \Tt^{\theta}$ satisfying this condition.
\epf

\begin{lemma}\label{lem:representatives}Let $w\in W^\theta$, $v\in W$,  $\w\in N_{\G^\theta}(\T)$ and $\v\in N_{\G}(\T)$ 
be representatives of $w$ and $v$, respectively. Let $y\in\G^\theta$ such that $y^{-1}F(y)=\w$. Then 
\begin{enumerate}
\item $\v \T\cap \G^\theta\cap \G^{F_w}\neq\emptyset$ if and only if $v\in 
W^\theta_w$.
\item An element $v$ in $W=N_{\Pad}(\Tad)/\Tad$ has a representative in $\Pad^{\theta,\circ}\cap \pi(\G^{F_w})$ if and 
only if $v\in W^\theta_w$. 
\end{enumerate}
\end{lemma}
\pf (1) If $\v \T\cap \G^\theta\neq\emptyset$, then, $\theta(\v)\in \v\T$, so $v\in W^\theta$ and we may assume 
$\v\in\G^\theta$. If $\v\T\cap\G^{F_w}\neq\emptyset$, then $F_w(\v)\in \v\T$, that is 
$\Ad(\w)(\v)\in\v\T$, {\it i.~e.} $wv=vw$. Conversely, assume $v\in W^\theta_w$. 
Now $W^\theta$ is the Weyl group of $\G^\theta$ and $F_w$ is a Steinberg endomorphism of 
$\G^\theta$ preserving its maximal torus $\T^\theta$.  By \cite[Proposition 23.2 
ff]{MT}, 
$$(W^\theta)^{F_w}=\left(N_{\G^\theta}(\T^\theta)/\T^\theta\right)^{F_w}\simeq 
N_{\G^\theta\cap\G^{F_w}}(\T^\theta)/(\T^\theta\cap\T^{F_w})$$
so any $v\in W_w^\theta=(W^\theta)^{F_w}$ has a representative in 
$$N_{\G^\theta\cap\G^{F_w}}(\T^\theta)=N_{\G^\theta}(\T^\theta)\cap \G^{F_w}=N_{\G^\theta}(\T)\cap 
\G^{F_w}=N_{\G}(\T)\cap \G^\theta\cap \G^{F_w}.$$

(2) Follows from (1) recalling that $\pi(\G^\theta)=\Pad^{\theta,\circ}$. 
\epf

We end the section with  a lemma that shows how some some of the results on the Weyl group apply  to the quest of preferred representatives in a twisted class. 

\begin{lemma}\label{lem:intersections} Let $t\in \Tt$ be such that $\Oc_t^{\theta,\G}\cap\PL_n(q)\neq\emptyset$. Then
\begin{enumerate}
\item There are $\sigma\in W^\theta$ and $\sigmad\in \sigma \Tt\cap \G^\theta$ such that  $\Oc_t^{\theta,\G}\cap 
\Tt_\sigma^F \neq\emptyset$. 
\item Let $\sigma$ be as in (1). Then $\Oc_t^{\theta,\G}\cap \Tt_w^F \neq\emptyset$ for every 
$w\in\Oc_{\sigma}^{W^\theta}$.
\item Fix $p$ odd and $x\in \Oc_t^{\theta,\G}\cap\PL_n(q)$. Then
$\Oc_t^{\theta,\G}\cap \PL_n(q)=\Oc_x^{\theta,\G^F}$.
\end{enumerate}
\end{lemma}
\pf(1) Pick a set of representatives $\{\dot{\tau},\;\tau\in W\}\subset N_{\G}(\Tt)$ as in Remark \ref{rem:representatives}. Let 
$g\in\G$ be such that $F(t)=gt\theta(g^{-1})$, see Remark \ref{rem:stable-orbit} (1). Let $u\in\U\cap \tau^{-1} 
\U^-\tau$, $\dot{\tau}\in N_{\G}(\Tt)\cap\tau \T$, $s\in\T$, 
$v\in \U$ such that $g=u\dot{\tau} sv$. Then 
\begin{align*}
F(t)\theta(g)=\left(F(t)\theta(u)F(t^{-1})\right)\cdot \left(
F(t)\theta(\dot{\tau})\theta(s)\right)\cdot \theta(v)\in 
\B\theta\tau\B. 
\end{align*}
On the other hand, $F(t)\theta(g)=gt=u\dot{\tau} svt=u (\dot{\tau} st )(t^{-1}vt)\in 
\B\tau\B$, which gives, by the uniqueness of the Bruhat decomposition, 
$\theta(\tau)=\tau\in W$ and, by construction, $\theta(\dot{\tau})=\dot{\tau}$. Also this yields 
$F(t)\theta(\dot{\tau})\theta(s)=\dot{\tau} st$, that is $F(t)=(\dot{\tau} 
s)\cdot_\theta t\in \Oc_t^{\theta,N_{\G}(\T)}$. Let $\sigmad:=\dot{\tau}^{-1}\in N_{\G^\theta}(\T)$. Then
$F_{\sigma}=\Ad(\sigmad)\circ F$ is again a Steinberg endomorphism for $\T$ and 
$F_{\sigma}(t)=ts\theta(s^{-1})\in \Oc_t^{\theta,\T}$. Let $r\in \T$ be such 
that $r^{-1}F_{\sigma}(r)=s$. Then $x=r^{-1}\cdot_\theta t\in 
\Oc_t^{\theta,\T}\cap \Tt^{F_{\sigma}}$. Indeed, 
\begin{align*}
 F_{\sigma}(x)=F_{\sigma}(r^{-1}) 
F_{\sigma}(t)\theta(F_{\sigma}(r))=F_{\sigma}(r^{-1}) 
st\theta(s^{-1}F_{\sigma}(r))=r^{-1}t\theta(r)=x.
\end{align*}
Let $y\in\G^\theta$ be such that $y^{-1}F(y)=\sigmad$ and set $z=y\cdot_\theta 
x=yxy^{-1}$. Then $z\in y\Tt^{F_\sigma}y^{-1}=(y\Tt y^{-1})^F\cap 
\Oc_t^{\theta,\G}$, by \eqref{eq:invariance} and (1) follows.

(2) By Lemma \ref{lem:invt} there is $g\in \SL_n(q)^\theta$ such that $g\Tt^F_\sigma 
g^{-1}=\Tt_w^F$. Hence, for $x\in \Oc_t^{\theta,\G}\cap \Tt_\sigma^F$ we have $g\cdot_\theta 
x\in  \Oc_t^{\theta,\G}\cap \Tt_w^F$.

(3) The group $\G_{\theta}(t)=\G^{\Ad(t^{-1})\circ\theta}$ is connected by \cite[Theorem 
8.1]{steinberg-endo} since $\Ad(t^{-1})\circ\theta$ is a semisimple 
automorphism as defined in \cite[p. 51]{steinberg-endo}. The result follows from \cite[Theorem 21.11]{MT}.
\epf

\section{Twisted classes of $\theta$-semisimple elements}\label{sec:regular}

We assume from now on that $p$ is odd. Recall the notation from \S \ref{sec:Gb}, specially in \eqref{eqn:gps}, \eqref{eqn:tori}.

\subsection{Strategy}

Next theorem is the first main result of the paper and a key step to apply the strategy in Section \ref{subsec:strategy}.
 
\begin{theorem}\label{thm:x-in-calT}
Let $x\in\PGL_n(q)$ be $\theta$-semisimple. Then there are $w\in 
W^\theta$ and $z\in \Tad^F_w$ such that $\Oc_x^{\theta,\Pb}=\Oc_z^{\theta,\Pb}$.
\end{theorem}
\pf
Let $\tx\in\PL_n(q)$ be such that $x=\pi(\tx)$. By Lemma \ref{lem:equivalence} (3), there is 
$g\in \G$ such that $g\cdot_\theta \tx=t\in \Tt$. Then there is $w\in W^\theta$ and 
$\tz\in \Oc_t^{\theta,\G}\cap \Tt_w^F$ such that $\Oc_\tx^{\theta,\G}=\Oc_\tz^{\theta,\G}$, 
by Lemma \ref{lem:intersections} (1). On the 
other hand, we have that 
$\Oc_\tz^{\theta,\G}\cap \PL_n(q)=\Oc_\tz^{\theta,\SL_n(q)}$, by Lemma \ref{lem:intersections} 
(3). The statement now follows applying $\pi$, for $z=\pi(\tz)$, as 
$\pi\left(\Tt^F_w\right)\subset \Tad^F_w$ and 
$\pi\left(\Oc_\tz^{\theta,\SL_n(q)}\right)=\Oc_z^{\theta,\Pb}$.
\epf

\subsubsection{The strategy}\label{subsec:strategy} 

Let $x$ be a $\theta$-semisimple element in $\PGL_n(q)$. By Theorem  \ref{thm:x-in-calT} we may assume $x\in\Tad_w^F$ for some $w\in W^\theta$. We have the 
following inclusions of subracks:
\begin{equation}\label{eqn:strategy-inclusions}
 \Oc_x^{\psi,\Pb}\supseteq \Oc_x^{\psi,\Pb}\cap\Tad_w^F\supseteq 
\Oc_x^{\psi,\pi(\T_w^F)}\simeq\Oc_1^{\psi,\pi(\T_w^F)}.
\end{equation}

We will establish sufficient conditions ensuring $\Oc_1^{\psi,\pi(\T_w^F)}$ is of type D. 
If the conditions are not satisfied and  
$\Oc_x^{\psi,\Pb}\cap\Tad_w^F\neq \Oc_x^{\psi,\pi(\T_w^F)}$, we will establish 
sufficient conditions ensuring $\Oc_x^{\psi,\Pb}$ is of type D. 

We look at the subracks 
$\Oc_x^{\theta,\pi(\T_w^F)}\simeq 
\Oc_1^{\theta,\pi(\T_w^F)}$ as in \eqref{eqn:strategy-inclusions}. Thus we investigate the 
abelian subgroups $\pi(\T_w^F)$ and $\pi(\T_w^F)\cap \Pb^\theta$. Let $\w\in w{\T}\cap 
\G^\theta$ and let $y\in \G^\theta$ be such that $y^{-1}F(y)=\w$.
We have 
$$\Oc_1^{\theta,\pi(\T_w^F)}\simeq 
\pi(\T_w^F)/(\pi(\T_w^F)\cap\Pb^\theta)\simeq \T_w^F/K$$
for
$K=\{t\in \T_w^F~|~\theta(t)\in t\Zc(\SL_n(q))\}$.
Let us set 
\begin{align*}
\tK_w=\{s\in 
\T^{F_w}~|~\theta(s)\in s\Zc(\SL_n(q))\}.
\end{align*}

\begin{lemma}\label{lem:zeta}
Assume $\Oc_x^{\theta,\Pb}\cap \Tad_w^F\neq\emptyset$. If there is $s\in \T^{F_w}/\tK_w$ such that $|s|$ is even and $>4$, 
then $\Oc_x^{\theta,\Pb}$ is of type D. 
\end{lemma}
\pf It follows from Remark \ref{rem:reductions}  (\ref{item:abelian2}) and Lemma \ref{lem:vendra}, as conjugation by 
$y$ gives the group isomorphism $\T^{F_w} /\tK_w\simeq \Oc_1^{\theta,\pi(\T_w^F)}$.
\epf

When conditions in Lemma \ref{lem:zeta} do not hold, we will 
use the following lemma.

\begin{lemma}\label{lem:two-orbits} Let $x\in\Tad_w^F$ for some $w\in 
W^\theta$, 
and assume $\Oc_x^{\theta,\Pb}\cap\Tad_w^F\neq\Oc_x^{\theta,\pi(\T_w^F)}$.
If there is  $z\in \mO_1^{\theta,\pi(\T_w^F)}\simeq\T^{F_w}/\tK_w$ such that $z^4\neq1$, then 
$\Oc_x^{\psi,\Pb}$ is of type D. 
\end{lemma}
\pf The subrack $X=\Oc_x^{\theta,\Pb}\cap\Tad_w^F$ is a disjoint union of orbits 
under the $\theta$-conjugation by $\pi(\T_w^F)$, one of which is 
$R=\Oc_x^{\theta,\pi(\T_w^F)}=x\Oc_1^{\theta,\pi(\T_w^F)}$. Let $S=\Oc_t^{\theta,\pi(\T_w^F)}=t\Oc_1^{\theta,\pi(\T_w^F)}\subset X$, $S\neq R$.
As $\Tad_w^F$ is abelian and $\theta^2=1$, \eqref{eq:rs} becomes
\begin{equation}\label{eqn:invo}(r\theta(r)^{-1})^2\neq(s\theta(s^{-1}))^2.\end{equation}
If  \eqref{eqn:invo} holds for $r:=x$, $s:=t$, we are done. Otherwise, we replace $s$ by $s'=sz\in S$, 
obtaining the desired inequality.\epf

\subsection{Conjugacy classes in $W^\theta$}\label{subsec:conj classes in W}

We need to describe $\tS_w$ and  $\tK_w$, $w\in W^\theta$.  We will use the 
identification of $W^\theta$ with $\sym_h\rtimes\Z_2^h$, for $h=\left[\frac{n}{2}\right]$. Set 
$\{\be_i:1\leq i\leq h\}$ the canonical $\Z_2$-basis of $\Z_2^h$. Also, for $\lambda=(\lambda_1,\ldots,\lambda_r)$ 
$\lambda_j\geq\lambda_{j+1}$ a partition of $h$, 
consider  the set ${\mathcal E}(\lambda)$ 
consisting of all vectors $\varepsilon\in \Z_2^r$ such that if 
$\lambda_j=\lambda_{j+1}$, then $\varepsilon_j=0$ implies 
$\varepsilon_{j+1}=0$. 

By Lemma 
\ref{lem:intersections} (2) it is enough to look at a set representatives of each 
$W^\theta$-conjugacy class. According to \cite[Proposition 24]{carter-cc} such 
a set is given by all 
$$
\sigma_{\lambda,\varepsilon}:=(1,2,\ldots,i_1)\be_{i_1}^{\varepsilon_1}(i_1+1, 
i_1+2,\ldots,i_2)\be_{i_2}^{\varepsilon_2}\cdots 
(i_{r-1},i_{r-1}+1,\ldots,h)\be_{h}^{\varepsilon_r}.
$$
with $i_j=\sum_{l\leq 
j}\lambda_j$ and $\varepsilon\in {\mathcal 
E}(\lambda)$.

To simplify the exposition, let $\vartheta:\I_n\to \I_n$ be the permutation 
$i\mapsto n+1-i$. Let us denote by $\bs_{p,q}$ the permutation $(p,q)$.
As an element in $\symm_{n}$, $w$ becomes a 
product of cycles as follows:
\begin{equation}\label{eqn:w-in-W}
w=\left(\bc_1\theta(\bc_1)\bs_{i_1,\vartheta(i_1)}^{\varepsilon_1}\right) \dots 
\left(\bc_h\theta(\bc_h)\bs_{i_h,\vartheta(i_h)}^{\varepsilon_h}\right),
\end{equation}
$\bc_j=(i_{j-1},i_{j-1}+1,\ldots,i_j)$, $1\leq j\leq h$, $i_{-1}=0$. We 
set $w_j:=\bc_j\theta(\bc_j)\bs_{i_j,\vartheta(i_j)}^{\varepsilon_j}$.

\smallbreak

We analyze cases $n$ odd and even separately and apply the results in Lemma 
\ref{lem:Tw}.

\subsubsection{$n$ odd}\label{subsubsec:n-odd}
Let $n=2h+1$ and $w=\sigma_{\lambda,\varepsilon}$.  
Let, for $j=1,\ldots, r$:
$$\F({j}):=\begin{cases}\F_{q^{\lambda_j}}^\times\times 
\F_{q^{\lambda_j}}^\times,& \text{ if 
}\varepsilon_j=0,\\
\F_{q^{2\lambda_j}}^\times, & \text{ if }\varepsilon_j=1.
\end{cases}$$
Direct computation shows  that $\Tt^{F_w}\simeq \F_q^\times\times 
\prod_{j=1}^r\F(j)$.
For $j\in \I_r$, $z_j\in \F(j)$, we set:
\begin{align*}
 &\overline{z}_j:=\begin{cases}x_jy_j,& \text{ if }\varepsilon_j=0 
\text{ and }z_j=(x_j,y_j),\\
z_j,& \text{ if }\varepsilon_j=1,
\end{cases}
&&\text{ and 
}&&\tz_j:=\overline{z}_j^{1+\varepsilon_jq^{\lambda_j}}\in\F_{q^{\lambda_j}}.
\end{align*}
 Observe that as $z_j$ 
runs in $\F(j)$ then 
$\tz_j$ covers $\F_{q^{\lambda_j}}^\times$ and $\tz_j^{(\lambda_j)_q}$ covers $\F_q^\times$. We have
$$\T^{F_w}:=\{(z,z_1,\ldots,z_r)\in \F_q^\times\times 
\prod_{j=1}^r\F(j)~|~z\prod_j\tz_j^{(\lambda_j)_q}=1\}\simeq 
\prod_{j=1}^r\F(j).$$
It follows from direct 
computation that 
\begin{align*}
\tK_w\simeq \{(z_1,\ldots,z_r)\in 
\prod_{j=1}^r\F(j)~|~\tz_j=\zeta, 1\leq j\leq r, \ 
\zeta\in \mu_n(\F_q)\}. 
\end{align*}
Hence, if $\gamma:\T^{F_w}\to \F_{q^{\lambda_1}}^\times \times \F_{q^{\lambda_1\lambda_2}}^\times 
\times \dots \times \F_{q^{\lambda_{r-1}\lambda_r}}^\times $ is given by
\begin{align}\label{eq:gamma}
 (z_1,\ldots,z_r)\mapsto (\tz_1^\bd, \tz_1\tz_2^{-1}, \dots,  
\tz_{r-1}\tz_r^{-1}),
\end{align}
then $\T^{F_w}/\tK_w\simeq \Imm \gamma$.

\subsubsection{$n$ even}\label{subsec:n-even}
Let $n=2h$, $w=\sigma_{\lambda,\varepsilon}$. With notation as in \S \ref{subsubsec:n-odd}, we have:
$$\T^{F_w}=\{(z_1,\ldots,z_r)\in
\prod_{j=1}^r\F(j)~|~\prod_j\tz_j^{(\lambda_j)_q}=1\}.$$
It follows from direct 
computation that 
\begin{align*}
\tK_w\simeq \{(z_1\ldots,z_r)\in\T^{F_w}~|~\tz_j=\zeta, 1\leq j\leq r, \ 
\zeta\in \mu_n(\F_q)\},
\end{align*}
hence $\T^{F_w}/\tK_w\simeq \Imm \gamma$, for $\gamma:\T^{F_w}\to \F_{q^{\lambda_1}}^\times \times 
\F_{q^{\lambda_1\lambda_2}}^\times 
\times \dots \times \F_{q^{\lambda_{r-1}\lambda_r}}^\times$ as in \eqref{eq:gamma}.

\subsection{Applying the strategy}We will deal with classes $\Oc_x^{\theta,\pi(\T_w^F)}$ for $x\in \Tad_w^F$. Observe 
that as $\Tad_w=\ty\Tad \ty^{-1}$, $x$ is represented by an element in 
$\Tt^{F_w}$ up to multiplication by matrices in $\Zc(\PL_n(q))$, {\it i.~e.} , up to a scalar factor 
in $\F_q^\times$. 
We apply Lemma \ref{lem:zeta} and the description of $\T^{F_w}/\tK_w$ from Section \ref{subsec:conj classes in W} on  each 
case to detect classes of type D. Let $\bone$ denote the partition 
$(1,\dots,1)$.

\begin{lemma}\label{lem:Tw}
Let  $\lambda=(\lambda_1,\dots,\lambda_r)$ be a partition of $h$, $\varepsilon\in{\mathcal E}(\lambda)$ and let $w=\sigma_{\lambda,\varepsilon}\in W^\theta$.
Let $x\in\Tad_w^F$. Then  
$\Oc_x^{\theta,\Pb}$ is of type D provided any of the following conditions hold.
\begin{enumerate}
 \item  $n$ is odd, $\lambda\neq \bone$.
\item $n$ is even, $\lambda\neq \bone$, and $r>2$.
\item $\lambda =\bone$, $n\neq3,4$ and $q>5$.
\item If $\lambda =\bone$, $n=3$ and $q=9,11$ or $q>13$. 
\item If $\lambda =\bone$, $n=4$ and $q>9$ and $q\equiv 1\mod(4)$.
\end{enumerate}
\end{lemma}
\pf In all cases we will  provide a suitable element in the image of the map $\gamma$ from \eqref{eq:gamma} and apply Lemma \ref{lem:zeta}.

(1) Assume $r>1$. If $j$ is such that $\varepsilon_j= 0$ 
and $\lambda_j>1$, consider $\tilde z_j=(x_j,1)$, for a generator $x_j\in 
\F_{q^{\lambda_j}}^\times$. If $$\gamma_j:=\gamma(1,\ldots,\tilde 
z_j,\ldots,1)=(x_j^{\delta_{1,j}\bd},\dots,x_j,x_j^{-1},\dots,1),$$ then 
$\mid\gamma_j\mid=\mid x_j\mid= q^{\lambda_j}-1>4$ and even. Similarly, if $r>1$ and $j$ is such 
that 
$\varepsilon_j= 1$ and $\lambda_j>1$, then it follows 
that if
$$\gamma_j:=\gamma(1,\ldots,z_j,\ldots,1)=
(z_j^{\delta_{1,j}\bd(1+q^{\lambda_j})},\dots,z_j^{1+q^{\lambda_j}},z_j^{-1-q^{\lambda_j}},\dots,
1)$$ for a generator $z_j$ of $\F_{q^{2\lambda_j}}^\times$, then 
$$\mid\gamma_j\mid=\mid 
z_j^{1+q^{\lambda_j}}\mid=\frac{q^{2\lambda_j}-1}{(q^{2\lambda_j}-1,1+q^{\lambda_j})}=q^{\lambda_j}-1>4.$$  
Now, if $r=1$, then $\lambda=(h)$, $h>1$. Pick $\overline{z}$ such that 
$\tz$ is a generator of $\F_{q^h}^\times$. Then 
$$
\mid\gamma(z)\mid=\mid\tz^\bd\mid=\frac{q^h-1}{(\bd,q^h-1)}=\frac{q-1}{\bd}(h)_q>(h)_q\geq 4.
$$
Observe that 
$\frac{q-1}{\bd}$ is always even,  whence the first inequality.  Moreover, $(h)_q=4$ only if $q=3$, $n=5$ in which case $\frac{q-1}{\bd}(h)_q=2(h)_q>4$.

(2) Assume now that $n$ is even.
We distinguish the following cases:

{\it Case $r>2$, $\lambda\neq\bone$.} 

Let us choose 
$\overline{z}_1$ such that $\tz_1$ is  a generator of $\F_{q^{\lambda_1}}^\times$. 
Choose $\tz_2=\dots=\tz_{r-1}=1$ and $z_r$ such that 
$\tz_1^{(\lambda_1)_q}\tz_r^{(\lambda_r)_q}=1$. Then $(z_1,\dots,z_r)\in\T^{F_w}$
and 
$$
\mid\gamma(z_1,\dots,z_r)\mid\geq\mid\tz_1\mid=q^{\lambda_1}-1>4 \text{ and even.}
$$ 

(3), (4), (5) If $n$ is odd, $n\neq3$ the computation in (1) shows that we 
can find $x\in\Imm\gamma$ with $\mid x\mid=q-1>4$ for $q>5$. If $n=3$, then $\Imm\gamma$ is cyclic of order $\frac{q-1}{\bd}>4$ for $q\geq9$, $q\neq 13$ and always even.

If $n$ is even, then $h=r\ge2$. If $r>2$ we may choose $\tz_1$ as a generator of 
$\F_{q}^\times$,
$\tz_2=\tz_1^{-1}$ and $\tz_j=1$ for $j\geq3$ and proceed as before. 
If $r=2$ then $n=4$. We need $\tz_2=\tz_1^{-1}$ and, choosing $\tz_1$ as above we have $|(\tz_1^\bd,\tz_1^2)|=\frac{q-1}{2}$. \epf

\begin{lemma}\label{lem:weyl}
Let $w\in W^\theta$ and $x\in \Tad_w^F$. If $\Oc_{w_0}^{W_w}$ is of type D, then 
$\Oc_{x}^{\theta,\Pb}$ is so. 
\end{lemma}
\pf Let $\w$ be a representative of $w$ in $\Pad^\theta\cap N(\Tad)$, see Remark \ref{rem:representatives}, and let $y\in (\Pad^\theta)^\circ=\pi(\G^\theta)$ be such that $y^{-1}F(y)=\w$, so $x=yty^{-1}$ for some 
$t\in \Tad^{F_w}$. Since $\Pb=[\Pad^F,\Pad^F]=y[\Pad^{F_w},\Pad^{F_w}]y^{-1}$ 
we have
$\Oc_{x}^{\theta,\Pb}\simeq \Oc_{t}^{\theta, [\Pad^{F_w},\Pad^{F_w}]}$. Now, 
${F_w}$ is again a Steinberg endomorphism of $\Pad$, 
and $\Tad$ is $F_w$-stable. Hence, \cite[Proposition 23.2]{MT} applies and by \cite[Exercise 
30.13]{MT} there is a group epimorphism 
$N_{\Pad}(\Tad)\cap [\Pad^{F_w},\Pad^{F_w}]\twoheadrightarrow W^{F_w}=W_w$ inducing a 
rack 
epimorphism  $\Oc_{x}^{\theta,\Pb}\twoheadrightarrow \Oc_1^{\theta,W_w}$. The statement
follows from Lemma \ref{lem:inner} (2).
\epf

For $\lambda=\bone$ and $j=0,\,\ldots,\,h$ we set 
$\varepsilon^j:=(\underbrace{1,\ldots,1}_\textrm{$(h-j)$ 
times},\underbrace{0,\ldots,0}_\textrm{$j$ times})\in {\mathcal E}(\lambda)$.

\begin{lemma}\label{lem:weyl-j}Let $w=\sigma_{\bone,\varepsilon^j}$ and let $x\in \Tad_w^F$. If $n$ is even and 
$j\geq3$, or if $n$ is odd and $j>3$, then $\Oc_{x}^{\theta,\Pb}$ is of type D. In particular, if  $x\in \Tad^F$, then  
$\Oc_{x}^{\theta,\Pb}$ is of type D provided $n\geq 6$, $n\neq 7$.
\end{lemma}
\pf By Lemma \ref{lem:weyl} it is enough to prove that $\Oc_{w_0}^{W_w}$ is of type D. 
Now $w\in W$ is the permutation 
$(1, n)\cdots(h-j, n+1-h+j)\in W'\times 1\leq W'\times W''$ where
$$W'\times W''\:=\sym_{\{1,\ldots,h-j,n+1-h+j,\ldots,n\}}\times\sym_{\{h-j+1,\ldots,n-h+j\}}\simeq \sym_{2(h-j)}\times \sym_{n-2(h-j)}$$ and 
$W_w=W'_w\times W''$, so
$\Oc_{w_0}^{W_w}\simeq \Oc_{w}^{W'_w}\times \Oc_{ww_0}^{W''}\simeq \Oc_{ww_0}^{W''}$. The latter  is of type D by \cite[Theorem 4.1]{AFGaV}.
\epf

\begin{lemma}\label{lem:blocks}
Assume $n=2h$ and let $\lambda=(\lambda_1,\dots,\lambda_h)$ be a partition 
of $h$. 
\begin{enumerate}
\item
If $w=\sigma_{\lambda,\varepsilon}=w_1\dots w_j\in W^\theta$ as in 
\eqref{eqn:w-in-W}, then there is a block matrix $y=\Diag(y_1,\dots,y_h)\in 
\G^{\theta}$ such that $\w=y^{-1}F(y)$ is a representative of $w$ in $N_{\G^\theta}(\T)\cap \SL_n(q)$, each block $y_j\in 
\Sp_{2\lambda_j}(\kk)$ and $\w_j=y_j^{-1}F(y_j)\in w_j \T$.
\item If $\lambda=(\lambda_1)$ and $w=\sigma_{\lambda,0}$, then there are $y_1\in 
\SL_{\lambda_1}(\kk)$ and $\w\in w\T\cap N_{\G^\theta}(\T)\cap \SL_n(q)$ such that $\w=y^{-1}F(y)$, $y=\Diag(y_1, \Jf_{\lambda_1}\,^t\!y_1^{-1}\Jf_{\lambda_1}^{-1})\in 
\G^{\theta}$.
\end{enumerate}
\end{lemma}
\pf
(1) Set $i_j=\sum_{l\leq 
j}\lambda_j$, $i_{-1}:=0$, $\Lambda_j=\{i_{j-1}+1,\dots,i_j\}$, $1\leq j\leq 
h$. Recall from \eqref{eqn:w-in-W} that $w\in \symm_{2h}$ can be viewed as an
element in $\symm_{2\lambda_1}\times \dots\times \symm_{2\lambda_h}$,
if we identify $\symm_{2\lambda_j}$ with 
the permutation group of $\Lambda_j\cup\vartheta(\Lambda_j)$, for $1\leq j\leq 
h$. Notice that $w_j=\bc_j\theta(\bc_j)\bs_{i_j,\vartheta(i_j)}\in 
\sym_{2\lambda_j}^{\theta_j}$ for each $1\leq j\leq h$. Hence each $w_j$ lies in the Weyl group of a $\theta$-invariant subgroup 
$\G_j\simeq \Sp_{2\lambda_j}(\kk)$ of $\G$, namely the subgroup of matrices of the 
shape
$$
\left(\begin{smallmatrix}
       \Id &     & & & \\
           & A& & B& \\
        & &\Id & & \\
     & C & & D& \\
        &     & & &\Id \\
      \end{smallmatrix}
\right), \qquad \left(\begin{smallmatrix}
       A & B \\
      C & D 
      \end{smallmatrix}\right)\in \Sp_{2\lambda_j}(\kk)
$$
and the non-zero entries outside the diagonal are indexed by integers in 
$\Lambda_j\cup\vartheta(\Lambda_j)$. Let us denote by $\theta_j$ the graph 
automorphism for $\G_j$. There exists 
a representative $\w_j$ of $w_j$ in $\G_j^{\theta_j}\simeq \Sp_{2\lambda_j}(\kk)$, 
as $n$ is even. Therefore, there exists $y_j\in \G_j^{\theta_j}\simeq 
\Sp_{2\lambda_j}(\kk)$ such that $y_j^{-1}F(y_j)=\w_j$. We remark that 
$[\G_i,\G_j]=1$ for $i\neq j$ and thus $y$ can be chosen as $y=y_1\dots y_h$.

(2) If $\varepsilon=0$ then $w$ lies in $\sym_{\lambda_1}$ and it is represented by block matrices of the form $\w=\Diag(A, \Jf_{\lambda_1}\,^t\!A^{-1}\Jf_{\lambda_1}^{-1})\in \G^{\theta}$. As we can always make sure that $A\in\SL_{\lambda_1}(q)$ \cite[Proposition 23.2]{MT}, we can apply Lang-Steinberg's Theorem to the connected group $\SL_{\lambda_1}(\kk)$.
\epf

\begin{lemma}\label{lem:n=2h,r=1,e=0}Let $n=2h$ for $h>1$, $\lambda=(h)$, $\varepsilon=(0)$ and 
 $w=\sigma_{\lambda,\varepsilon}$. Let $x\in \Tad_w^F$. Then $\Oc_x^{\theta,\Pb}$ is of 
type D provided one of the following holds:
\begin{enumerate}
\item $x\theta$ is not an involution and $\Pb\neq \PSL_4(3), \PSL_4(7)$.
\item $n\geq 6$.
\item $x\theta$ is an involution, $n=4$ and $q\equiv 1(4)$, $q\neq5,\,9$.
\end{enumerate} 
\end{lemma}

\pf 
We have 
$w=(1,2,\cdots,h)(n,n-1,\cdots,h+1)$. Let $y\in \Sp_n(\kk)$ satisfy
$\w=\pi(y^{-1}F(y))$. Set $\ty=\pi(y)$. Thus we may assume 
$$x=\ty\pi(t)\ty^{-1}, \quad \text{ for }  
t=\diag(a,a^q,\dots,a^{q^{h-1}},b^{q^{h-1}},\dots,b),$$
for some $a,b\in\F_{q^h}^\times$. We set, for $\xi\in \kk, 
\xi^{(h)_q}=1$:
$$t_\xi=\diag(a\xi,(a\xi)^q,\dots,(a\xi)^{q^{h-1}},(b\xi)^{q^{h-1}},\dots,
b\xi)\in \Tad_w^F.$$ It follows that 
$\Oc^{\theta,\pi(\T_w^F)}_x=\ty\{\pi(t_\xi):\xi\in \kk, 
\xi^{(h)_q}=1\}\ty^{-1}$.

Set $\dag=\dag_n:=\pi(\diag(-\id_h,\id_h))\in\PGL_n(q)$. Notice that $x\theta$ 
is an involution if and only if $\theta x=x^{-1}$ which happens only if
$x\in 
\Oc_1^{\theta,\Tad_w^F}\cup \Oc_\dag^{\theta,\Tad_w^F}$. 
We claim that if $x\not\in 
\Oc_1^{\theta,\Tad_w^F}\cup \Oc_\dag^{\theta,\Tad_w^F}$, then 
$\Oc_x^{\theta,\Pb}\cap \Tad_w^F\neq\Oc_x^{\theta,\pi(\T_w^F)}$.

Let us compute the (twisted) action of $w_0\in W_w^\theta=\langle w,w_0\rangle$ on 
$\pi(t)$.
We have
\begin{align*}
 w_0\cdot_\theta 
\pi(t)t&=\diag(b,b^q,\dots,b^{q^{h-1}},a^{q^{h-1}},\dots,a).
\end{align*}
Hence, 
$w_0\cdot_\theta x\in 
\Oc^{\theta,\pi(\T_w^F)}_x$ only if $ab^{-1}=ba^{-1}$. This gives the claim.

\medbreak
(1)
We apply 
Lemma \ref{lem:two-orbits}: we search for $z\in \mO_1^{\theta,\pi(\T_w^F)}$ 
such 
that $\mid z\mid\neq 1,2,4$.
According to the discussion in \S \ref{subsec:n-even}, 
$\mO_1^{\theta,\pi(\T_w^F)}$ is a cyclic group of order $\ell$, for 
$\ell=\frac{(h)_q}{(\bd,(h)_q)}=\frac{(h)_q}{(q-1,h)}$, as $\bd=(q-1,2h)$ and 
$(q-1,(h)_q)=(q-1,h)$. If $h=2$, so $n=4$, we have $\ell=\frac{1+q}{2}$, so $q\neq 3,7$ is enough.
If $h$ is odd, then $\ell$ is odd and $\ell>1$ since 
$\ell>\frac{1+q}{q-1}$. Then we can find such a $z$. From now we shall assume that $h\geq 4$ is even. 
We distinguish three cases, according to 
$h>q-1$, $h=q-1$ or $h<q-1$. 
If $h>q-1$ then
$
\ell>\frac{1+q(h-1)}{q-1}=h+\frac{h-(q-1)}{q-1}>4
$ 
and we are done. The same computation proves the claim if $h=q-1>4$. If $h=q-1=4$ a direct computation gives the claim. 
Finally, if $h<q-1$, then
$\ell>\frac{(h)_{h+1}}{h}\geq \frac{h+h(h-2)_h}{h}>6$. 

(2) If $x\theta$ is not an involution, then we apply (1). If $x\theta$ is an involution, then we
have that either $\Oc_x^{\theta,\Pb}\simeq \Oc_1^{\theta,\Pb}$ or 
$\Oc_x^{\theta,\Pb}\simeq 
\Oc_{\dag}^{\theta,\Pb}$, by Lemma \ref{lem:inner} (3). Now, $1,\kappa\in \Tad^F$  by Lemma \ref{lem:blocks} (2), as we may assume  $y=\Diag(A,\Jf \,^tA^{-1}\Jf^{-1})$ for some there $A\in\PL_h(\kk)$.  
Thus, by Lemma \ref{lem:weyl-j}, $\Oc_x^{\theta,\Pb}$ is of 
type D if $h\geq 3$.

(3) Since $1,\kappa\in \Tad^F$, we apply Lemma \ref{lem:Tw} (5).
\epf

\begin{proposition}\label{pro:h2}
Let $q\equiv 3(4)$, $q\neq 3,7$, $\Pb=PSL_4(q)$. Let $t$ be either $1$ or 
$\kappa=\left(\begin{smallmatrix}
         \id_2&0\\0&-\id_2
        \end{smallmatrix}\right)$. Then $\Oc_t^{\theta,\Pb}$
is of type D.
\end{proposition}
\pf We will apply Lemma~\ref{lem:vendra}. It is enough to find  $x\in \Pb$ such 
that the order of  $xt\theta(x)^{-1}t$ in $\Pb$ is even and $>4$.
Set ${\bf u}_t:\Pb\to \Pb$, ${\bf u}_t(x)=-xt\theta(x)^{-1}t=xt\Jf \,^t x \Jf 
t$. 
For each $e,f\in\F_q$, $A,E,F\in\F_q^{2\times 2}$, with $E,F$ traceless, let us set 
\begin{align*}
&{\bf m}(A,e,f)=\left(\begin{matrix}A&e\,\id_2\\
f\,\id_2&\Jf_2 \,^tA\Jf_2^{-1}
\end{matrix}\right),&&
{\bf n}(A,E,F)=\left(\begin{matrix}A&E\\
F&\Jf_2 \,^tA\Jf_2^{-1}\end{matrix}\right).
\end{align*}
We have that ${\bf u}_1({\bf m}(A,e,f))={\bf m}(A,e,f)^2$, ${\bf u}_\kappa({\bf 
n}(A,E,F))={\bf n}(A,E,F)^2$.

Moreover, for any $x\in\Pb$ there are $e,f\in\F_q$, $A,E,F\in\F_q^{2\times 2}$, 
$E,F$ traceless such that ${\bf u}_1(x)={\bf m}(A,e,f)$ 
and ${\bf u}_\kappa(x)={\bf m}(A,E,F)$.

We shall exhibit a matrix ${\bf m}(A,0,0)={\bf n}(A,0,0)$ whose projective order 
is a multiple of $4$ and it is bigger than $8$. This will prove the statement 
for both $t=1,\kappa$.

Let $\F^\times_{q^2}=\langle \xi \rangle$ and consider 
the matrix 
$z=\diag(\xi^{\frac{q-1}{2}},-\xi^{\frac{1-q}{2}},-\xi^{\frac{1-q}{2}},\xi^{
\frac{q-1}{2}}
)$ in $\SL_4(\F_{q^2})$. The order of $z$ is $2({q+1})$ and 
$z^{\frac{q+1}{2}}=\diag(\omega,\omega^{-1},\omega^{-1},\omega)$
for $\omega$ a primitive fourth root of $1$, hence the projective order of $z$ 
is $q+1$. 

We claim that  $z$ is $\PGL_4(\overline{\Fq})$-conjugate to 
$x={\bf m}(\left(\begin{smallmatrix}\Tr(z)/2&1\\
1&0\\
\end{smallmatrix}\right),0,0)$ and that $\Tr(z)\in\F_q^\times$.  If this is the 
case, 
${\bf u}_1(x)=x^2$ and its projective order is $\frac{q+1}{2}$ which is
even as $q\equiv 3(4)$ and bigger than 4 since $q\geq11$. 

The claim  is proved if the
following conditions hold, namely
\begin{align*}
&\det z=1; &&
\Tr z=2(\xi^{\frac{q-1}{2}}-\xi^{\frac{1-q}{2}})\in\F_q;
&& \xi^{\frac{q-1}2}\neq -\xi^{\frac{1-q}2}.
\end{align*}
Indeed, in this case, the matrix ${\bf m}(\left(\begin{smallmatrix}\Tr(z)/2&1\\
1&0\\
\end{smallmatrix}\right),0,0)$ is diagonalizable and it is necessarily 
$\PL_4(\overline{\Fq})$-conjugate to the matrix $z$.
The first and third conditions are immediate. For the second, let $\sigma$ be 
the
(involutive) generator of the Galois group $\Gal(\F_{q^2},\F_q)$ of the  
extension
$\F_q\subset \F_{q^2}$. We need $\sigma(\Tr x)=\Tr x$. But $\sigma$ coincides 
with
$\Fr_p^m$, that is $\sigma(\xi)=\xi^q$ and thus the equality holds. 
\epf

\begin{lemma}\label{lem:n=2h,r=1,e=1}Let $n=2h$, $h>1$,  $\lambda=(h)$ and 
$\varepsilon=(1)$, 
$w=\sigma_{\lambda,\varepsilon}$. Let $x\in \Tad_w^F$. Then $\Oc_x^{\theta,\Pb}$ is of 
type D provided that one of the following holds:
\begin{enumerate}
\item $x\theta$ is not an involution and $\Pb\neq \PSL_4(3), \PSL_4(7)$.
\item $x\in \Oc_1^{\theta,\PGL_n(q)}$, $n\geq6$.
\item $x\in \Oc_1^{\theta,\PGL_n(q)}$, $n=4$, $q>9$.
\item $x\theta$ is an involution, $x\not\in\Oc_1^{\theta,\PGL_n(q)}$ and $h$ is even. 
\end{enumerate}
\end{lemma}
\pf  In this case $w=(1,2,\ldots,h-1,h,n,n-1,\ldots,h+2,h+1)$ as a permutation in $\sym_n$. Arguing as in Lemma \ref{lem:n=2h,r=1,e=0}
we may assume that, for some $a\in\F_{q^{n}}^\times$,
$$x=\ty\pi(t)\ty^{-1}, \quad \text{ for }  
t=t_a=\diag(a,a^q,\dots,a^{q^{h-1}},a^{q^{2h-1}},\dots,a^{q^h})$$ and $\ty$ such that $\ty^{-1}F(\ty)=\w$.

(1) Notice that $x\theta$ is an involution if and only if $a^2$ lies in 
$\F_{q^h}^\times$. Now, set, for $\xi\in \F_{q^n}^\times$ such that $\xi^{(n)_q}=1$ and $z=\xi^{1+q^h}$ in $C_{(h)_q}\subset \F_{q^h}^\times$:
$$t_{az}=\diag(az,(az)^q,\dots,(az)^{q^{h-1}},\dots,
(az)^{q^h})\in \Tt^{F_w}.$$ It follows that 
$\Oc^{\theta,\pi(\T_w^F)}_x=\ty\{\pi(t_{az}):\xi\in \kk, 
\xi^{(n)_q}=1\}\ty^{-1}$.  
Observe that $(yw_0^{-1}y^{-1})\cdot_\theta x=y \pi(t^{q^h}) y^{-1}$ lies in 
$\Oc_x^{\theta,\pi(\T_w^F)}$ if and only if $a^2\in \F_{q^h}^\times$. In other 
words, if $w_0\cdot_\theta t$ lies in $\Oc_t^{\theta, \pi(\T^{F_w})}$ only if 
$x\theta$ is an involution.
If this is not the case, we can proceed as in Lemma \ref{lem:n=2h,r=1,e=0} and obtain that if $\Pb\neq 
\PSL_4(3), \PSL_4(7)$, then $\Oc_x^{\theta,\Pb}$ is of type D.

(2), (3) Assume $a\in\F_{q^h}^\times$. Then $\Oc_x^{\theta,\Tad_w^F}=\Oc_1^{\theta,\Tad_w^F}$ and 
$\Oc_x^{\theta,\pi(\T_w^F)}\cap\Tad_w^F=\Oc_x^{\theta,\pi(\T_w^F)}$. 
In this case $x\in\Oc_1^{\theta, (\ty \Tad\ty^{-1})^F}\subset 
\Oc_1^{\theta,\PGL_n(q)}$ and we may assume $a=1$, $t=\id$ by Proposition \ref{pro:x-in-PGL}. If $n\geq6$ this class is of type D by Lemma \ref{lem:weyl-j}. Assume $n=4$. If $q\equiv1(4)$, $q>9$, then we may apply Lemma \ref{lem:Tw} (5). For $q\equiv 3(4)$, $q>7$ we apply Proposition \ref{pro:h2}

(4) Assume $a^{q^h}=-a$ and moreover that $h$ is even. We apply 
Lemma \ref{lem:vendra}: We search for $r\in\Oc_{x\theta}^{\Pb\ltimes\langle\theta\rangle}$ such that 
$\mid rx\theta\mid$ is even and bigger than 4. Equivalently, we look for $z\in \Pb^{F_w}$ such that 
the order of $(z\cdot_{\theta} t)\theta \, t\theta=(z\cdot_{\theta} t) t^{-1}$ is even and bigger than 4. If $h$ is 
even, then this is achieved by taking 
$z=t$, 
as $t\cdot_{\theta} t=t^3$ and thus $\mid t^2\mid=(h)_q$ which is even and bigger than 4, as $h>1$.
\epf

\bigskip

We are missing the
case $r=2$, $\lambda=(\lambda_1,\lambda_2)\neq \bone$. That is, the case in 
which 
$n=2(\lambda_1+\lambda_2)$, with $\lambda_2\geq 1$,  $\lambda_1\geq 2$. This 
is the content of Lemmas \ref{lem:r=2,varepsilon=0} (when $\varepsilon_1=0$) and 
\ref{lem:r=2,varepsilon=1} (when $\varepsilon_1=1$). Let us set 
$$w=c_1\theta(c_1)s_{\lambda_1,n-\lambda_1+1}^{\varepsilon_1}c_2\theta(c_2)s_{
\lambda_1+\lambda_2,\lambda_1+\lambda_2+1}^{\varepsilon_2}$$ where
$c_1=(1,2,\ldots,\lambda_1)$ and $c_2=(\lambda_1+1,\ldots,\lambda_1+\lambda_2)$. 
We write 
$$w_1:=c_1\theta(c_1)s_{\lambda_1,n-\lambda_1+1}^{\varepsilon_1},\quad
w_2=c_2\theta(c_2)s_{\lambda_1+\lambda_2,\lambda_1+\lambda_2+1}^{\varepsilon_2}
.$$ The group $W^\theta_w$ always contains $w_1$, $w_2$ and the elements
$$w_{0}^{(1)}=(1n)\cdots (\lambda_1,n-\lambda_1+1),\quad 
w_{0}^{(2)}=(\lambda_1+1,n-\lambda_1)\cdots 
(\lambda_1+\lambda_2,\lambda_1+\lambda_2+1),$$ which correspond to the longest 
elements in each block.

\begin{lemma}\label{lem:r=2,varepsilon=0}
Let $n=2(\lambda_1+\lambda_2)$, $\lambda=(\lambda_1,\lambda_2)$ with $\lambda_1\geq2$, $\lambda_2\geq1$, 
and $\varepsilon_1=0$. Let $x\in 
\Tad_w^F$. If  $q>5$, then
$\Oc_x^{\theta,\Pb}$ is of 
type D.
\end{lemma}
\pf We have $x=\pi(\ty)\pi(t)\pi(\ty)^{-1}$ where $\ty^{-1}F(\ty)=\w$,  
$t=\diag(t_{11},t_2,t_{12})$. Here $t_{11}$ and $t_{12}$ are diagonal matrices 
of size $\lambda_1$ and $t_2\in \PL_{2\lambda_2}(\ku)^{F_{w_2}}$ is also 
diagonal. Also, $t_1:=\diag(t_{11},t_{12})\in 
\PL_{2\lambda_1}(\ku)^{F_{w_1}}$.

Assume first $W_w^\theta\cdot \Oc_{\pi(t)}^{\theta,\pi(\T^{F_w})}\neq 
\Oc_{\pi(t)}^{\theta,\pi(\T^{F_w})}$. 
With notation as in \S \ref{subsec:n-even}, 
the abelian group
$\Oc_1^{\theta,\pi(\T_w^F)}$ is isomorphic to the image 
of the map $\gamma$. Since $\lambda_1>1$ we may choose $z_1$ so that 
$|\gamma(z_1,1)|=|(\tz_1^{\bf d}, \tz_1)|=(\lambda_1)_q>4$ and Lemma 
\ref{lem:two-orbits} applies. 

We now determine  when $W_w^\theta\cdot \Oc_{\pi(t)}^{\theta,\pi(\T^{F_w})}= 
\Oc_{\pi(t)}^{\theta,\pi(\T^{F_w})}$, acting by $w_0^{(1)}$. Arguing as in 
Lemma \ref{lem:n=2h,r=1,e=0}, we see that 
$w_0^{(1)}\cdot_\theta \pi(t)\in\Oc_{\pi(t)}^{\theta,\pi(\T^{F_w})}$ only if $\Oc_x^{\theta,\Pb}\simeq \Oc_{x'}^{\theta,\Pb}$ where $x'$ has the form
$x'=\pi(\ty\diag(\id_{\lambda_1},t_2',\pm \id_{\lambda_1})\ty^{-1})$ with 
$t_2'$ a  diagonal element in $\PL_{2\lambda_2}(\kk)^{F_{w_2}}$. By Lemma 
\ref{lem:blocks} (2), $x'$ lies in $\Tad_{w_2}^F$. So, if $\lambda_2=1$, then the partition associated with $w_2$ is $\lambda'=\bone$ and Lemma \ref{lem:Tw} (3) applies. 
If $\lambda_2>1$, then the partition associated with $w_2$ is 
$\lambda'\neq\bone$ so $r=\lambda_1+1>2$ and  Lemma \ref{lem:Tw} (2) applies.
\epf

\begin{remark}\label{rem:r=2}
It follows from the proof of Lemma \ref{lem:r=2,varepsilon=0} that even in the 
case 
$q=5$ the class
$\Oc_x^{\theta,\Pb}$ is of type D, provided $\lambda_2>1$. Also, if $q=3$ then 
this class is
of type D provided $\lambda_1>1$ and $\lambda_2>2$. 
\end{remark}

\begin{lemma}\label{lem:r=2,varepsilon=1}
Let $n=2(\lambda_1+\lambda_2)$, $\lambda=(\lambda_1,\lambda_2)$ with $\lambda_1\geq2$, $\lambda_2\geq1$, 
and $\varepsilon_1=1$. Let $x\in 
\Tad_w^F$. Then
$\Oc_x^{\theta,\Pb}$ is of 
type D.
\end{lemma}
\pf We follow the strategy and notation from Lemma \ref{lem:r=2,varepsilon=0}. 
In this case  $x=\pi(\ty)\pi(t)\pi(\ty)^{-1}$ for 
$t=\diag(t_{11},t_2,t_{12})$, and 
$$t_1=\diag(t_{11},t_{12})=(a,a^q,\ldots,a^{q^{\lambda_1-1}},a^{q^{2\lambda_1-1}
} , \ldots , a^{q^{\lambda_1}})\in\PL_{2\lambda_1}(\kk)^{F_{w_1}}.$$
Applying $w_0^{(1)}$ and arguing as in the proof of Lemma \ref{lem:n=2h,r=1,e=1}, we see that $\w_0^{(1)}\cdot_\theta \Oc_{\pi(t)}^{\theta,\pi(\T^{F_{w}})}\neq \Oc_{\pi(t)}^{\theta,\pi(\T^{F_{w}})}$ with the possible exception of the case in which $a^2\in\F_{q^\lambda_1}^\times$.
If $a\in\F^\times_{q^{\lambda_1}}$, then there are diagonal elements $d_{11}$, 
$d_{12}$ in $\PL_{\lambda_1}(\kk)$ and $d_2\in 
\PL_{2\lambda_2}(\kk)^{F_{w_2}}$ such that 
$d_1:=\diag(d_{11},d_{12})\in \PL_{2\lambda_1}(\kk)^{F_{w_1}}$, 
$\det d=1$ for $d:=\diag(d_{11},d_2,d_{12})$ and 
$d\cdot_\theta t=(\id_{\lambda_1},t_2',\id_{\lambda_1})$. 
By Lemma \ref{lem:blocks} the latter lies in the $F$-stable maximal torus 
associated with $w_2$. 
In this case, $r>2$ and we apply Lemma \ref{lem:Tw} (2). We assume from now on 
that $a^2\in\F_{q^\lambda_1}^\times$ and $a\not\in\F_{q^\lambda_1}^\times$, {\it 
i.~e.}, $a^{q^{\lambda_1}-1}=-1$. 
Hence if $\varepsilon_2=0$, then, the element $t$ is:
\begin{align*}
t=t_{a,b,c}=(a,a^q,\ldots,a^{q^{\lambda_1-1}},b,b^q,\ldots,b^{q^{\lambda_2-1}},
c^ { q^{\lambda_2-1}},\ldots,c, -a^{q^{\lambda_1-1}},\ldots,-a),
\end{align*}
 for $b,c\in\F_{q^{\lambda_2}}^\times$, while if 
$\varepsilon_2=1$, then 
\begin{align*}
t=t_{a,b}=(a,a^q,\ldots,a^{q^{\lambda_1-1}},b,\ldots,b^{q^{\lambda_2-1}},b^{
q^ {2\lambda_2-1}},\ldots,b^{q^\lambda_2}, 
-a^{q^{\lambda_1-1}},\ldots,-a),
\end{align*}
for $b\in \F_{q^{2\lambda_2}}^\times$.
The $\T^{F_w}$-orbit consists of elements of the form $t_{a\tz_1,b\tz_2,c\tz_2}$ ($t_{a\tz_1,b\tz_2}$, respectively) for $\tz_1\in \F_{q^{\lambda_1}}^\times$ and $\tz_2\in\F_{q^{\lambda_2}}^\times$ satisfying $\tz_1^{(\lambda_1)_q}\tz_2^{(\lambda_2)_q}=1$. 
Since all elements in $\F_{q^{2\lambda_1}}^\times$ satisfying $a^{q^{\lambda_1}-1}=-1$ lie in the same $\F_{q^{\lambda_1}}^\times$-coset, we may assume that $|a|=2(q^{\lambda_1}-1)$. We consider $w^{-1}_1\cdot \pi(t_{a,b})$ ($w^{-1}_1\cdot \pi(t_{a,b, c})$, respectively). If it lies in a different $\pi(\T^{F_w})$-orbit 
than $\pi(t_{a,b})$  ($\pi(t_{a,b, c})$, respectively),  we apply Lemma \ref{lem:two-orbits}. Otherwise, there are
$\tz_1\in \F_{q^{\lambda_1}}^\times$, $\ell\in\F_q^\times$, and $\tz_2\in \F_{q^{\lambda_2}}^\times$ such that $\ell\tz_2=1$, $a^{q-1}=\ell \tz_1$, and $\tz_1^{(\lambda_1)_q}\tz_2^{(\lambda_2)_q}=1$. If this is the case, then 
$|a^{q-1}|=|\tz_2^{-1} \tz_1|=2(\lambda_1)_q$. Thus, there is an element in the image of the map $\gamma$ in \eqref{eq:gamma} of even order $>4$ and Lemma \ref{lem:zeta} applies.\epf

\begin{remark}\label{rem:missing}
Lemma \ref{lem:n=2h,r=1,e=1} does not cover the case $h$ odd, $W_w^\theta\cdot_\theta \Oc_{\pi(t)}^{\theta,\pi(\T^{F_w})}=\Oc_{\pi(t)}^{\theta,\pi(\T^{F_w})}$, and $x\neq 1$. This actually amounts to at most a single class for each group, up to rack isomorphism: Keep the notation from the lemma, 
let $\zeta$ be a generator of $\F_{q^n}^\times$, $\eta=\zeta^{\frac{1+q^h}{2}}$ and set
\begin{align}\label{eqn:missing}
\eth=\ty \pi(t_\eta)\ty^{-1}.
\end{align}
Then $x\in \Oc_\eth^{\theta, \PGL_n(q)}$ and $\Oc_x^{\theta,\Pb}\simeq \Oc_\eth^{\theta,\Pb}$, by Proposition \ref{pro:x-in-PGL}. 

(1) This class seems to be difficult. For instance, the subrack 
$$\Oc_{\eth}^{\theta, \pi(\T_w^F)}\simeq \Oc_{\pi(t_\eta)}^{\theta, \pi(\T^{F_w})}= \Oc_{\pi(t_\eta)}^{\theta, N_{\pi(\G^{F_w,\theta})}(\T)}$$ is not of type D. 
Indeed, $|\pi(t_\eta) \theta|=2$  and  $\Oc_{t_\eta}^{\theta, \pi(\Tt^{F_w})}=\{\pi(t_{az}):\xi\in \kk, 
\xi^{(n)_q}=1,z=\xi^{1+q^h}\}$. Then for every $r\in \Oc_{t_\eta}^{\theta, \pi(\Tt^{F_w})}$ the order  
$|r\theta \pi(t_\eta)\theta|=|r\pi(t_\eta^{-1})|$ divides  $(h)_q$ and hence it is odd. 

(2) This class does not occur when dealing with $\theta$-conjugacy classes 
in $\Pb$ instead of $\PGL_n(q)$, if $q\equiv1(4)$. Let $s=t_a$ as in the proof of Lemma \ref{lem:n=2h,r=1,e=1}. Then $s$ lies in $\SL_n(\kk)$ only if $q\equiv 
3(4)$. Indeed $\det(s)=a^{(n)_q}=a^{(1+q^h)(h)_q}=-a^{2(h)_q}=1$ gives  $a^{2(h)_q}=-1$.  Also, $-1=a^{(1-q)(h)_q}=a^{\frac{q-1}{2} 
2(h)_q}=(-1)^{\frac{q-1}{2}}$ so $\frac{q-1}{2}$ is odd.  

\end{remark}

In group-theoretical terms, the class of $\eth$ is of type D if and only if the following question has an affirmative 
answer:
\begin{question}\label{question}
Let $\eta$ and $s=t_\eta$ be as above, recall the matrix $\Jf$ from \eqref{eqn:J}. 
Is there a matrix in $A\in\SL_n(\overline{\F_q})^{F_w}$ such that 
the projective order of 
\begin{align}\label{eqn:m(A)}
 m(A):=\Jf A\Jf  \,^t\!A s
\end{align}
is even and bigger than 4?
\end{question}

\subsubsection{A class not of of type D}

So far, we have seen that most $\theta$-semisimple classes in $\PSL_n(\Fq)$ are of type D. 
Next proposition shows that there exist classes that are not of type D. This shows that the condition $q\neq3$ in 
Proposition \ref{pro:h2} is necessary.

\begin{proposition}\label{pro:not}
The class $\Oc_1^{\theta,\PSL_4(3)}$
is not of type D.
\end{proposition}
\pf
We adopt the notation from the proof of Proposition \ref{pro:h2} to show 
that the projective order of any matrix $Y$ of the form ${\bf m}(A,e,f)$ is at 
most equal to $4$.
For such $Y$, we verify that $\det Y=(\det A-ef)^2$, that $\Tr(Y)=2\Tr(A)$, 
and that the matrix $Y$ annihilates the polynomial $X^2-\Tr(A)X+(\det A-ef)$. 
Since $Y\in{\rm SL}_4(\F_3)$ we have $\delta:=\det A-ef=\pm1$. Therefore, whenever $\Tr(A)=0$, 
the matrix $Y$ is an involution in $\PSL_4(3)$. Let us assume $\Tr(A)\neq0$. 
We have $Y^2=\Tr(A)Y-\delta$ so 
$$Y^3=\Tr(A)Y^2-\delta Y=\Tr(A)^2 Y-\delta\Tr(A)-\delta 
Y=(1-\delta)Y-\delta\Tr(A).$$
If $\delta=1$, then $Y^3=\pm1$. If $\delta=-1$, then
$Y^4=-Y^2+\Tr(A)Y=-1$.\epf

\subsection{Proof of Theorem \ref{thm:one}}\label{the-proof}

\pf
We cite the rows in the table according to their position in the last 
column. We make the convention that the head row is the 0th one. Let us consider the class of an element $x$
that might possibly be not of type D.
If $\lambda\neq\bone$ then by Lemma \ref{lem:Tw} (1) and (2), $n$ has to be even and $r\leq 2$. 
By Lemmata \ref{lem:r=2,varepsilon=0} and \ref{lem:r=2,varepsilon=1}, the case $r=2$ and either $\varepsilon_1=1$
or $q>5$ is ruled out, yielding the first row.
If $r=1$ and $\varepsilon=(0)$, then by Lemma \ref{lem:n=2h,r=1,e=0} and Proposition \ref{pro:h2},
the classes that are not of type D may occur only for $n=4$ and  either 
$q=3,7$, or $\theta(x)=x^{-1}$ and $q=5,9$. This gives the second and the third row. 
If $r=1$ and $\varepsilon=(1)$, then by Lemma \ref{lem:n=2h,r=1,e=1}, 
the classes that are not of type D may occur only in the following situations:
\begin{itemize} 
\item $n=4$; $q=3,7$ and $\theta(x)\neq x^{-1}$, which is the fourth row;
\item $n=4$, $q=3,5, 7,9$ and  $x\in\Oc_1^{\theta,\PGL_n(q)}$. This case is considered in the last row, as, 
up to rack isomorphism, this class is 
represented by an element in an $F$-stable maximal torus with associated partition $\bone$. 
\item  $n$ twice an odd number, $x\in\Oc_\nu^{\theta,\PGL_n(q)}$, which is the fifth row. 
\end{itemize}
Assume now $\lambda=\bone$. Then, by Lemma \ref{lem:Tw} the only classes that could occur in the table are those for
$q=3,5$, or for $n=3$ and 
$q=7,13$, or else $n=4$ and $q\equiv 3(4)$, or $q=9$. This gives the sixth and the seventh row.
\epf

\section{Twisted classes of elements with trivial $\theta$-semisimple part}\label{sec:unipotent}

Recall that for $x\in\PGL_n(q)$ there is a unique decomposition $x=us$ with $u$ unipotent, $s$ a 
$\theta$-semisimple element and 
$us=s\theta(u)$. Then  we have the rack inclusions, see \eqref{eqn:unipotent-inclusions}:
$$\Oc_x^{\theta,\Pb}\supset\Oc_u^{\Pb_{\theta}(s)}.$$

Assume $s=\pi(s')$, $s'\in \Tt$.
According to \cite[8.1]{steinberg-endo}, $\G_{\theta}(s')$ is a connected reductive group. By \cite[Theorem 1.1, Proposition 3.1]{mohr-tesis},  any simple component of $\G_{\theta}(s')$ is isomorphic either to  $Sp_{2a}(\ku)$ or to $\SL_a(\ku)$, $a\in\N$, if $n$ is even; and either to $Sp_{2a}(\ku)$,  $\SL_a(\ku)$ or  $SO_{2a+1}(\ku)$, $a\in \N$, if $n$ is odd. Taking $F$-invariants and arguing as in \cite[\S 3]{carter-centralizers},  one sees that $\G_{\theta}(s')^F$ contains a product of finite classical groups: unitary, special linear and symplectic if $n$ is even, and 
orthogonal, unitary, special linear and symplectic if $n$ is odd. Then $\Oc_x^{\theta,\Pb}$ is of type D whenever the conjugacy class of some component of  $u$ in one of these factors is so. 

The group $\Pb_{\theta}(s)=\pi(\G^F)^{Ad(s)\circ \theta}$ might properly contain $\pi((\G^{Ad(s')\circ \theta})^F)=\pi(\G_{\theta}(s')^F)$ although the latter already contains all unipotent elements. So, even if $\Oc_u^{\theta, \pi(\G_{\theta}(s')^F)}$ fails to be of type D,  it is still possible that $\Oc_u^{\theta,\Pb_{\theta}(s)}$ is so.

The unipotent classes in $\PSL_n(q)$ and $Sp_{2n}(q)$ are studied in \cite{ACGa, ACGa2}, whereas a similar analysis for unitary groups and orthogonal groups is in preparation.  This enables us to draw conclusions in case $s=1$. 

\begin{proposition}\label{pro:trivial-even}
Assume $n=2h$ is even and $q>3$. 
Let $\Oc$ be a $\theta$-twisted conjugacy class with trivial $\theta$-semisimple part and non-trivial unipotent part. Then $\Oc$ is of type D. 
\end{proposition}
\pf By the discussion in Subsection \ref{subsec:semis}, a representative of the class is a unipotent element $u\in\Pb^\theta$ and $\Oc$ has a 
subrack isomorphic to $\Oc_u^{\id,(\Pb^\theta)^\circ}\simeq\Oc_u^{\id,Sp_n(q)}$ (see the isogeny argument in \cite[1.2]{ACGa}). This rack is of type D with 
the exception of the classes with Jordan form corresponding to the partition 
$(2,1,\ldots,1)$, for $q$ either $9$ or not a square.
The reader should be alert that the form used for defining $Sp_n(\kk)$ in 
\cite{ACGa2} differs from the one considered 
here. Explicitly, they are 
related by the change of basis:
$$
e_i\mapsto\begin{cases}e_i& i\le h \mbox{ odd } \mbox { or } i> h \mbox{ even;} \\
e_{n-i+1}& i\le h \mbox{ even } \mbox{ or } i> h \mbox{ odd.}\\
\end{cases}
$$
Hence, if $M$ is the matrix that gives this basis change, elements in $\G^\theta$ are obtained from matrices therein 
conjugating by $M$. 

Let us consider the partition $(2,1,\ldots,1)$. There are two unipotent classes in $\Sp_n(q)$ associated with it. They are represented by 
$u_1=1+\eta_1e_{1,n}$ and by $u_2=1+\eta_2 e_{1,n}$ for $\eta_1$  a square and $\eta_2$ not a square in $\F_q^\times$, respectively. 
Consider $g=\pi(\diag(-\id_{2},\id_{n-2}))\in\Pb$ and set 
$$v_i:=g\cdot_\theta 
u_i=\pi\left(\sum_{\begin{smallmatrix}j=1,2,\\ \,n-1,n\end{smallmatrix}}e_{jj}-\sum_{j=3}^{n-3}e_{jj}-\eta_i e_{1n}\right)\in \Pb^\theta, \qquad i=1,2.$$ 
Let us consider the matrix $\sigma=\left(\begin{smallmatrix} 0&&1\\
&\id_{n-2}\\
-1&&0\\
\end{smallmatrix}\right)$; in particular $\pi(\sigma)\in \Pb^\theta$. 

Now, set $r_i=u_i$, $s_i=\pi(\sigma)\cdot_\theta v_i$ and $R_i=\Oc_{u_i}^{\id,\Pb^\theta}$, 
$S_i=\Oc_{s_i}^{\id,\Pb^\theta}$, $i=1,2$. It follows that the elements $r_i, s_i$ satisfy 
\eqref{eqn:typeD}. Unless $q\equiv 1(4)$ and $n=4$, the subracks $R$ and $S$ are disjoint as one is a unipotent class and the other is not. 

Assume now $n=4$, $q\equiv1(4)$. 
Then, for $i=1,2$ there is $\xi_i\in \F_q^\times$ such that $\xi_i^2\neq1$ and 
$(\xi\eta_i)^2\neq2$. Indeed, this excludes at most 4 elements, hence the case 
of $q>5$ follows, whereas if $q=5$, then 2 is not a square and we can take 
$\xi_i=2\in\F_5$. 
In this case we take $g=\pi(\diag(1,\xi_i,1,\xi_i^{-1}))\in\Pb$ and 
$$
v_i:=g\cdot_\theta u_i=\pi\left(\begin{smallmatrix}
1&&\xi_i\eta_i\\
&\xi_i^2\id_2\\
&&1\end{smallmatrix}\right), \qquad i=1,2.$$
Then we consider $r=u_i$, $s=\pi(s)\cdot_\theta v_i$ with $s_i$ as above, $R=\Oc_{u_i}^{\id,\Pb^\theta}$ and 
$S=\Oc_{s_i}^{\id,\Pb^\theta}$ and the proposition follows.
\epf

If $n$ is odd, less is known about unipotent conjugacy classes in $\Pb^\theta$. We can still obtain the following. Recall that in this case $\G^\theta\simeq\SO_n(\ku)$.
\begin{proposition}\label{pro:trivial-odd}
Assume $n>3$ is odd. 
Let $\Oc$ be a $\theta$-twisted conjugacy class with trivial $\theta$-semisimple part. If the Jordan form of its $p$-part in $\Pb^\theta$ corresponds to the partition $(n)$, then $\Oc$ is of type D.
\end{proposition}
\pf As above, a representative of the class is $u\in\Pb^\theta$ and $\Oc$ has a subrack isomorphic to $\Oc_u^{\id,\Pb^\theta}\simeq\Oc_u^{\id,SO_n(q)}$. We apply \cite[3.7]{ACGa2}.
\epf

\section*{Acknowledgements}

The authors are grateful to  Nicol\'as Andruskiewitsch, Kei-Yuen Chan, Fernando Fantino, Gast\'on Garc\'ia, Donna Testerman and Leandro
Vendramin for useful discussions and important e-mail exchanges. We are also indebted to Ernest Vinberg for the suggestion in Remark \ref{rem:vinberg} and to the referee for useful comments on the exposition of this manuscript.

\end{document}